\font\bb=msbm10

\def\bN{\hbox{\bb N}}
\def\bP{\hbox{\bb P}}
\def\bQ{\hbox{\bb Q}}
\def\bR{\hbox{\bb R}}

\def\bZ{\hbox{\bb Z}}
\def\b1{{\bf 1}}

\def\cB{{\cal B}}
\def\cD{{\cal D}}
\def\cF{{\cal F}}
\def\cG{{\cal G}}

\def\cI{{\cal I}}
\def\cJ{{\cal J}}
\def\cL{{\cal L}}
\def\cN{{\cal N}}
\def\cZ{{\cal Z}}

\def\endproof{\nobreak$\qquad$\nobreak$\sqcap$\llap{$\sqcup$}}

\def\leaderfill{\leaders\hbox to 1em{\hss.\hss}\hfill}

\centerline{\bf LOCAL FIELDS, GAUSSIAN MEASURES, AND BROWNIAN MOTIONS}
\vskip 24pt
\centerline{\it Steven N. Evans\footnote*{Preparation supported in part by a Presidential Young Investigator Award and Alfred P. Sloan Foundation Fellowship}}
\centerline{\it University of California at Berkeley}
\centerline{\it Department of Statistics}
\centerline{\it 367 Evans Hall}
\centerline{\it Berkeley, CA 94720-3860}
\centerline{\it USA}
\centerline{\it (evans@stat.berkeley.edu)}

\vskip 30pt

\leftline{1. Introduction}
\vskip 18pt
\leftline{2. Local fields}
\vskip 18pt
\leftline{3. Normed spaces and orthogonality}
\vskip 18pt
\leftline{4. Properties of $K$-Gaussian probability measures}
\vskip 18pt
\leftline{5. Construction and first properties of $K$-Brownian motion}
\vskip 18pt
\leftline{6. Random measures, polar sets, and $K$-Brownian motion}
\vskip 18pt
\leftline{7. Constructions and continuity of $K$-Brownian local time}
\vskip 18pt
\leftline{8. Other $K$-Gaussian random series}
\vskip 18pt
\leftline{~~~Appendix}
\vskip 18pt
\leftline{~~~References}
\vfill\eject

\centerline{\bf 1. Introduction}

\vskip 24pt

 A local field is any locally compact, non-discrete field other than the field of real numbers or the field of complex numbers.  All local fields are totally disconnected, and are either finite algebraic extensions of the field of $p$-adic numbers or finite algebraic extensions of the the less familiar $p$-series field (the field of formal Laurent series with coefficients drawn from the finite field with $p$ elements.).  

   Local fields and the vector spaces over them have a rich algebraic and geometric structure that in some ways rivals that of the corresponding objects for the real or complex fields.  These notes are a survey of the author's attempts to find and understand canonical probabilistic entities in a local field setting.  In particular, we propose answers to the related questions,\hfil\break
\centerline{``What are the analogues for Gaussian measures?''}
and \hfil\break
\centerline{``What are the analogues for Brownian motion and its multiparameter relatives?''}
We have discussed these problems in [Evans, 1989a, 1991, 1993], and this work is a distillation  of those papers.

      The original motivation for the study of local fields comes from number theory (cf. [Cassels, 1986]) and the theory of group representations (cf. Chap. XII of [Curtis and Reiner, 1962]).  However, local fields have found a multitude of other applications, from the the study of error-free computation (cf. [Gregory and Krishnamurthy, 1984]) to mathematical physics (cf. [Khrennikov, 1990], [Ruelle and Thiran, 1989], [Spokoiny, 1989], [Vladimirov and Volovich, 1989] and [Brydges et al., 1992] and the references therein).

   Local fields have also become  objects of interest in their own right.  In particular, there is a well-developed theory of analysis in this context.   A representative sampling of the area may be found in [Dwork, 1982], [Iwasawa, 1972], [Koblitz, 1980], [Mahler, 1980], [Monna, 1970], [Schikhof, 1984], [Taibleson, 1975] and [van Rooij, 1978]. We will draw on some elements of this work, but our goal is to make our treatment as self-contained as possible by at least quoting in full most of the results we need.

   In contrast, there has not been a similarly extensive study of probability on local field objects.  Most of the small body of work that we are aware of may be found in [Albeverio and Karwowski, 1991, 1994], 
[Brillinger, 1991], [Evans, 1988a, 1988b, 1989a, 1989b, 1991, 1993], [Guimier, 1989], [Madrecki, 1983, 1985, 1990, 1991], and [Missarov, 1989, 1991].  We should remark, however, that if one ignores the algebraic structure of local fields and thinks of them merely as ultrametric spaces or sets with a tree-like structure, then this work can be seen as part of the large and growing literature on probability in such a setting.   Also, we note that consideration of processes with an ultrametric index set has played a fundamental role in the study of boundedness and continuity of general Gaussian processes (see, for example, Chap. IV of [Adler, 1990]). 

   It soon becomes apparent in the study of (Euclidean) Gaussian measures that they are the
class of probability measures which is forced upon us if we search for a
probabilistic structure that is tightly linked to the linearity and orthogonality
properties of the spaces on which the measures are defined.
For instance, one entry point into the theory of Gaussian random variables on
an arbitrary real vector space with suitable measurable structure is to
define a random variable $X$ as being Gaussian if whenever $X_1$, $X_2$ are
two independent copies of $X$, then the pair
$( \alpha_{11} X_1 + \alpha_{12} X_2 , \alpha_{21} X_1 + \alpha_{22} X_2 )$
has the same law as $(X_1 , X_2 )$ for each pair of orthonormal vectors
$( \alpha_{11} , \alpha_{12} )$, $( \alpha_{21} , \alpha_{22} ) \in \bR^2$.
It can be shown that, in the appropriate special cases, this abstract definition
is equivalent to the usual concrete definitions for $\bR^n$-valued
Gaussian random variables and Gaussian stochastic processes.

  There is a suitable concept of orthogonality in the local field setting,
and so we can mimic the abstract definition given above and see where it takes us.  We remark that essentially the same class of measures that we are led to was derived in [Madrecki, 1983] from different considerations.  (As far as we can ascertain, this paper has not been published, but some of its essence can be gleaned from [Madrecki, 1985, 1990, 1991].)

  The resulting theory is similar in many ways to the Euclidean theory.
For example, linear transformations of ``Gaussian'' variables are ``Gaussian'' and
orthogonality is equivalent to independence.
We also observe the sort of zero-one behaviour that we have come to expect
from the Gaussian theory, absolute continuity questions are easily understood, and there is an analogue of the body of ideas around  Cameron-Martin type theorems and the related concept of reproducing kernel Hilbert space.

   Unfortunately, we also lose something in our new setting.
Roughly speaking, the local field notion of orthogonality is an $\cL^{\infty}$
object rather than an $\cL^2$ object.
As a consequence, the resulting ``Gaussian'' theory is not a second-order
theory where some analogue for the concept of covariance describes the distributional
picture.
In particular, the law of a ``Gaussian'' process is not described by its family
of two-dimensional marginal distributions.

   Ordinary Brownian motion and its multiparameter relatives, such as L\'evy's multiparameter Brownian motion and the Brownian sheet, are Gaussian processes with  covariance structures that are intimately tied to the algebraic and geometric structure of the underlying Euclidean index set.  Our goal of answering the second question above is thus to construct and study an analogous local field Gaussian process that takes values in a local field vector space, is indexed by another local field vector space, and is defined in such a way that there is a similar close linkage between the metric and dependence properties of the process and the properties of the index set.  The analogy between this {\it local field Brownian motion} and ordinary Brownian motion is not, of course, complete.  For example, there is no natural order structure on a local field, and so many of the features of the Euclidean theory in which order plays a prominent role (for example, the martingale property or the analytic theory of transition semigroups and resolvents) do not appear here.

     Our hope is that this new process will turn out to be as useful in the investigation of local field objects  as ordinary Brownian motion is in  Euclidean settings as diverse as classical potential theory, partial differential equations and complex analysis. We make one small step in that direction by establishing an intimate connection between 
the characterisation of polar sets for local field Brownian motion 
and the theory of Riesz potentials on a local field that has previously been investigated in [Taibleson, 1975].  This connection is similar to
the connection between polar sets for Brownian motion
and Newtonian or logarithmic potential theory discovered
 by [Kakutani, 1944a, 1944b].

   Our efforts to characterise the polar sets lead us to study the
 notion of additive functionals or homogeneous random measures.  We show that, as in the Brownian case, each such random measure is uniquely described by a (deterministic) `characteristic measure' on the state space.
This measure is the counterpart of the Revuz measure in the
Brownian case. We give necessary and sufficient conditions 
for a given measure to be a characteristic measure.

   We study the counterpart of the homogeneous random measure {\it par excellence} - Brownian local time.  When it exists, the local time at a point in the state space can be recovered  by intrinsic constructions on the corresponding level set that parallel, respectively, the dilation construction of [Kingman, 1973] and the Hausdorff measure construction of [Taylor and Wendel, 1966].  Moreover, there is an analogue of the 
theorem  of [Trotter, 1958]; namely, the local time is continuous in ``space and time''.

   The local field Brownian motion can be constructed as
the sum a series of deterministic functions with independent,
local field Gaussian 
coefficients.  We finish by saying a little 
about other such random series. 
Although we show that there are broad classes of random series
which are stationary, we also find that there is no 
obvious counterpart
to the representation of a general stationary Gaussian process on the circle
as a random Fourier series.

\vfill\eject

\centerline{\bf 2. Local fields}

\vskip 24pt

This section is essentially a summary of selected results from [Taibleson, 1975]
and [Schikhof, 1984].
We refer the reader to these works for a fuller account.  Before giving the
general definition of a local field, we begin with the prototypical
example.

\vskip 24pt
\noindent
{\bf Example.}
Fix a positive prime $p$.
We can write any non-zero rational number $r \in \bQ \backslash \{ 0 \}$ uniquely as $r=p^s (a/b)$
where $a$ and $b$ are not divisible by $p$. Set $\vert r \vert = p^{-s}$.
If we set $\vert 0 \vert =0$, then the map $\vert \cdot \vert$ has
the properties:
$$
\leqalignno{&\vert x \vert = 0 \Leftrightarrow x=0;&(2.1)\cr
&\vert xy \vert = \vert x \vert \vert y \vert ;&(2.2)\cr
&\vert x+y \vert \le \vert x \vert \vee \vert y \vert .&(2.3)\cr}
$$
The map $(x,y) \mapsto \vert x-y \vert$ defines a metric on $\bQ$ and we
denote the completion of $\bQ$ in this metric by $\bQ_p$.
The field operations on $\bQ$ extend continuously  to make $\bQ_p$
a  topological field called {\it $p$-adic numbers}. 
The map $\vert \cdot \vert$ also extends continuously and the extension
has properties (2.1) - (2.3).  The closed unit ball around $0$ 
$\bZ_p = \{x \in \bQ_p : |x| \le 1 \}$ 
is the closure in $\bQ_p$ of the integers $\bZ$, 
and is thus
a ring (this is also apparent from (2.1) - (2.3)),
called the  {\it $p$-adic integers}.
As
$\bZ_p = \{x \in \bQ_p : |x| < p\}$, the set $\bZ_p$ is also open.
Any other ball around $0$ is of the form
$\{x \in \bQ_p : |x| \le p^{-k}\} = p^k \bZ_p$
for some integer $k$.  Such a ball is the closure of the rational numbers
divisible by $p^k$, and is thus a $\bZ_p$ - module (this is again also apparent from (2.1) - (2.3)).  In particular, such a ball is an additive
subgroup of $\bQ_p$.  Arbitrary
balls are translates (=  cosets) of these closed and open subgroups.
In particular, the topology of $\bQ_p$ has a base of closed and open sets,
and hence $\bQ_p$ is totally disconnected.  
Further, each of these balls is compact, and hence $\bQ_p$
is also locally compact.

\vskip 24pt

   A {\it local field} is a locally compact, non-discrete, totally disconnected, topological field.  (As an aside, a locally compact, non-discrete, topological field that
is not totally disconnected is necessarily 
either the real or the complex numbers.  Also, every local field is either a finite algebraic extension of the $p$-adic number field for some prime $p$ or a finite algebraic extension of the $p$-{\it series field}; that is, the field of formal Laurent series with coefficients drawn from the finite field with $p$ elements.)

  From now on, we let $K$ be a fixed local field.
There is a real-valued mapping on $K$ which we denote by $x \mapsto \vert x \vert$.
This map has the properties (2.1) - (2.3) and it takes the values 
$\{ q^k : k \in \bZ \}
\cup \{ 0 \}$, where $q = p^c$ for some prime $p$ and positive integer $c$
(so that for $K=\bQ_p$ we have $c=1$).

A map with properties (2.1)-(2.3) is called a non-archimedean valuation.
Property (2.3) is known as the {\it ultrametric inequality} or
the {\it strong triangle inequality}.
The mapping $(x,y) \mapsto \vert x-y \vert$ on $K \times K$ is a metric on $K$
which gives the topology of $K$.
A consequence of (2.3) is that if $\vert x \vert \ne \vert y \vert$, then $\vert x+y \vert = \vert x \vert \vee \vert y \vert$.
This latter result implies that for every ``triangle''
$\{ x,y,z \} \subset K$ we have that at least two of the lengths $\vert x-y \vert$,
$\vert x-z \vert$, $\vert y-z \vert$ must be equal and is therefore often
called the {\it isosceles triangle property}.

  In the words of [Schikhof, 1984], ``... we shall follow a bad but widespread
habit and omit the subscript $K$ in $1_K$ and $n_K$ (:= the sum of $n$ times $1_K$).''
Clearly $\vert 1 \vert =1$.
We will write $D$ for $\{ x: \vert x \vert \le 1 \}$ (so that $D=\bZ_p$
when $K=\bQ_p$).
If we choose $\rho \in K$ so that $\vert \rho \vert =q^{-1}$, then
$$
\eqalign{\rho^k D & = \{ x: \vert x \vert
\le q^{-k} \} \cr
& = \{ x: \vert x \vert < q^{-(k-1)} \} \cr}
$$
for each $k \in \bZ$ (so that for $K=\bQ_p$ we could take $\rho = p$).

 The set $D$ is a ring  (the so-called {\it ring of integers} of $K$). Each of the sets $\rho^k D$, $k \in \bZ$, is a compact $D$- submodule of $K$ and every
non-trivial compact $D$-submodule of $K$ is of this form. For $\ell <k$ the
additive quotient group $\rho^\ell D / \rho^k D$ has order $q^{k-\ell}$.  Consequently, $D$ is
the union of $q$ disjoint translates of $\rho D$.  Each of these
components is, in turn, the union of $q$ disjoint translates
of $\rho^2 D$, and so on.  We can thus think of the collection of
balls contained in $D$ as being arranged in an infinite rooted $q$-ary tree:
the root is $D$ itself, the nodes at level $k$ are the balls of radius
$q^{-k}$ (= cosets of $\rho^k D$), and the $q$ ``children'' of such a ball
are the $q$ cosets of $\rho^{k+1} D$ that it contains.  We can uniquely
 associate
each point in $D$ with the sequence of balls that contain it, and so we
can think of the points in $D$ as the boundary of this tree.  This tree
picture alone does not capture all the algebraic structure of $D$; the rings of
integers for the $p$-adic
numbers and the $p$-series field are both represented by a $p$-ary tree,
even though the $p$-adic field has characteristic $0$ whereas the $p$-series
field has characteristic $p$. 

  Analogous tree structures are inherent in certain ``self-similar'' Euclidean fractals.  For example, the points of the classical Cantor set 
and the Sierpinski 
gasket can be  identified naturally with the boundary of the 
infinite binary and ternary
trees respectively.  The resulting depiction of local fields as
Euclidean fractals is discussed in [Cuoco, 1991].

  We may make this ``parameterisation" of $D$ as the boundary of the infinite, rooted, $q$-ary
tree a little more explicit by relating it to
 the usual identification of the boundary of
the infinite, rooted, $q$-ary tree with the set of infinite sequences
drawn from a set of $q$ elements.
Let $\{a_1, \ldots, a_q\}$ be a
complete list of coset representatives of $\rho D$
in $D$ (so that for $K=\bQ_p$ we could take $\{0,1,2, \ldots, p-1\}$).  We may express each
$y \in D$ uniquely as $y= \sum_{j=0}^{\infty} b_j \rho^j$, where
$b_j \in \{ a_1 , \ldots , a_q \}$ for each $j$.

   There is a unique Borel measure $\lambda$
on $K$ for which
$$
\eqalign{\lambda (x+A) = \lambda (A) ,~~~~~~~~~~ & x \in K , \cr
\lambda (xA) = \vert x \vert \lambda (A),~~~~~~~~~~ & x \in K , \cr}
$$
and
$$
\lambda (D) =1 ,
$$
The measure $\lambda$ is a suitably normalised Haar measure on the additive group of $K$.   For ease of notation, we will often write $dx$ for $\lambda(dx)$.  In the case of $\bQ_p$, the restriction of $\lambda$ to $\bZ_p$ is the weak limit as $n \rightarrow \infty$ of the sequence of probability measures that at the $n$-th stage assigns mass $p^{-n}$  to each of the points  $\{0,1,\ldots,p^n-1\}$.

   We will take our normalisation of Haar measure on the additive group of $K^n$
to be such that we have
the product of $n$-copies of $\lambda$.  With a slight abuse of
notation, we will also denote this measure by $\lambda$ if the context
is clear, and also write $dx$ for $\lambda(dx)$. 

There is a character $\chi$ on the additive group of $K$ with the properties
$$
\chi( D ) = \{ 1 \} ,
$$
and
$$
\chi(  \rho^{-1} D ) \ne \{ 1 \} .
$$
For $n=1,2, \ldots$, the correspondence $\xi \mapsto \chi_{\xi}$,
where $\chi_{\xi} (x) = \chi( \xi \cdot x )$, establishes an isomorphism
between the additive group of $K^n$ and its dual.
The uniqueness theorem for Fourier transforms in this setting thus states that
if $\nu_1$, $\nu_2$ are two finite measures on $K^n$ such that
$$
\int \chi( \xi \cdot x) \nu_1 (dx) = \int \chi( \xi \cdot x ) \nu_2 (dx)
$$
for all $\xi \in K^n$, then $\nu_1 = \nu_2$.

There is one Fourier transform which is of particular interest to us.
If $\Phi : [0, \infty [ \rightarrow \{ 0,1 \}$ is the indicator function
of the interval $[0,1]$, then
$$
q^{-n} \int_K \chi( \xi x) \Phi (q^{-n} \vert x \vert ) \lambda (dx) = \Phi (q^n \vert \xi \vert ) . \leqno(2.4)
$$
We remark that $\Phi$ has the property
$$
\Phi (a \vee b)= \Phi (a) \Phi (b), ~~~~~a,b \ge 0 . \leqno(2.5)
$$

\vfill\eject

\vskip 30pt
\noindent
\centerline{\bf 3.  Normed spaces and orthogonality}

\vskip 24pt

The material in this section is included for ease of reference and is
a summary of results and ideas which may be found in [Schikhof, 1984].

\vskip 24pt
\noindent
{\bf Definition 3.1.}
Let $E$ be a vector space over $K$.
A norm on $E$ is a map $\Vert~~~\Vert : E \rightarrow [ 0, \infty [$ such that
$$
\leqalignno{&\Vert x \Vert = 0 \Leftrightarrow x=0;&(3.2)\cr
&\Vert \lambda x \Vert = \vert \lambda \vert \Vert x \Vert ,~~\lambda \in K ,&(3.3)\cr
&\Vert x+y \Vert \le \Vert x \Vert \vee \Vert y \Vert .&(3.4)\cr}
$$
We call the pair $(E, \Vert~~~\Vert )$ a normed vector space (over $K$).
If $E$ is complete in the metric $(x,y) \mapsto \Vert x-y \Vert$, we
say that $E$ is a Banach space (over $K$).

Property (3.4) is also called the ultrametric inequality and leads to the
obvious analogue of the isosceles triangle property.

\vskip 24pt
\noindent
{\bf Example 3.5.}
All normed finite-dimensional vector spaces over $K$ are Banach spaces.
In particular, for $n=1,2, \ldots ,$ the space $(K^n , \vert~~~\vert )$, where
$$
\vert (x_1 , \ldots , x_n ) \vert = \vert x_1 \vert \vee \ldots \vee \vert x_n \vert ,
$$
is a Banach space.  The balls around $0$ in this space are all of the
form
$$\{(x_1, \ldots, x_n) : |(x_1, \ldots, x_n)| \le q^{-k}\}
= \{(x_1, \ldots, x_n) : |x_i| \le q^{-k}, \; 1 \le i \le k\}
= (\rho^k D)^N
= \rho^k D^N
$$
for some integer $k$.

\vskip 24pt
\noindent
{\bf Example 3.6.}
Let $T$ be a compact topological space.
The space $C(T,K)$ of $K$-valued continuous functions on $T$ equipped with
the norm $\Vert~~~\Vert_C$ given by
$$
\Vert f \Vert_C = \sup \{ \vert f(t) \vert : t \in T \}
$$
is a Banach space.  Example 3.5 is just the special case when $T$
has $n$ points.

\vskip 24pt
\noindent
{\bf Example 3.7.}
Let $( \Omega , {\cal F} , \bP)$ be a probability space.
Let $\cL^{\infty}(\bP)$ be the set of measurable functions $f: \Omega \rightarrow K$
such that ${\rm ess} \sup \{ \vert f( \omega ) \vert : \omega \in \Omega \} < \infty$.
If we say that $f=g$ whenever $f( \omega ) = g( \omega )$ for almost all $\omega$,
then $\cL^{\infty}(\bP)$ equipped with the norm $\Vert~~~\Vert_{\infty}$ defined by
$$
\Vert f \Vert_{\infty} = {\rm ess} \sup \{ \vert f( \omega ) \vert : \omega \in \Omega \}
$$
is a Banach space.  Example 3.5 is the special case when $\Omega$ has
$n$ points and $\bP$ is any measure that assigns positive mass to every
point in $\Omega$.
(This example is not in [Schikhof, 1984], but the proof proceeds just as in the
parallel real case.)

\vskip 24pt

    We now want to introduce a possible definition for orthogonality in our local field setting.  In order to motivate this definition,
we recall that it is possible
to characterise orthogonal subsets of a real inner product space
using only the associated norm and not explicitly making use of the inner product.  Let $(H,\langle \; \rangle)$ be a real inner product
space with associated norm $\| \; \|_H$, then a susbset $G \subset H$ will
be orthogonal if and only if for every finite subset $\{y_1, \ldots, y_n\} \subset G$ and each collection of scalars $\beta_1, \ldots \beta_n \in \bR$
we have
$$\|\sum_{i=1}^n \beta_i y_i \|_H = (\sum_{i=1}^n |\beta_i|^2 \|y_i\|_H^2)^{1/2}.$$
The quantity on the right-hand side is just the Euclidean norm of the vector
$(|\beta_1| \|y_1\|_H, \ldots, |\beta_n| \|y_n\|_H) \in \bR^n$.  The
following definition mimics this characterisation by replacing the Euclidean norm with the ``canonical'' norm on $K^n$ introduced in Example 3.5.  This
is the first instance of a phenomenon that will repeatedly recoccur: the
the most fruitful analogues of real or complex constructs in the local field
setting are often obtained by replacing ``$\cL^2$'' by ``$\cL^\infty$''

\vskip 24pt
\noindent
{\bf Definition 3.8.}
Suppose that $(E, \Vert~~~\Vert )$ is a normed space (over $K$).
We say that a set $F \subset E$ is orthogonal if for every finite subset
$\{ x_1 , \ldots x_n \} \subset F$ and each $\alpha_1 , \ldots , \alpha_n \in K$, we have
$$
\Vert \sum_{i=1}^n \alpha_i x_i \Vert = \bigvee_{i=1}^n \vert \alpha_i \vert \Vert x_i \Vert .
$$
We say that an orthogonal set $F \subset E$ is orthonormal if $\Vert x \Vert =1$
for all $x \in F$.

\vskip 24pt
   If $E$ is a separable Banach space then $E$ has a countable orthogonal basis.  This basis may be taken to be orthonormal if $\{\|x\| : x \in E\} = \{|\alpha|: \alpha \in K\} = \{q^k\}_{k \in \bZ} \cup \{0\}$.

\vskip 24pt
\noindent
{\bf Example 3.9.}
For $x \in \bZ_p$ and $n \in \{ 1,2, \ldots \}$, set
$$
{x \choose n} = {x(x-1) \ldots (x-n+1) \over n!} .
$$
Set ${x \choose 0} =1$.
Then the functions ${\cdot \choose 0}$, ${\cdot \choose 1}$, ${\cdot \choose 2} , \ldots$
form an orthonormal basis (the Mahler basis) for $(C( \bZ_p, \bQ_p )$, $\Vert~~~\Vert_C )$

\vskip 24pt
\noindent
{\bf Example 3.10.}
Recall from \S 2 that we can write $x \in \bZ_p$ as $x= \sum_{k=0}^{\infty} b_k p^k$,
where $b_k \in \{ 0,1, \ldots , p-1 \}$ for all $k$.
Given $m \in \{ 1,2, \ldots \}$, we write $m \triangleleft x$ if
$m= \sum_{k=0}^N b_kp^k$ for some $N \in \bN$ and adopt the convention
that $0 \triangleleft x$ for all $x$.
If $n \in \{ 1,2, \ldots \}$, then
$$
\{ m : m \triangleleft n , m \ne n \}
$$
is finite and has a largest element (in the order defined by the relation $\triangleleft$)
which we denote by $n_{\rm \_}$.
The functions $e_0 , e_1 , \ldots$ defined by
$$
e_n (x) = \cases{1,&if $n \triangleleft x$ ,\cr
0,&otherwise,\cr}
$$
form an orthonormal basis (the van der Put basis) for $C( \bZ_p, \bQ_p )$.
If $f \in C( \bZ_p, \bQ_p )$ has the expansion
$$
f(x) = \sum_{n=0}^{\infty} a_n e_n (x) ,
$$
then $a_0 =f(0)$ and $a_n =f(n)-f(n_{\rm \_} )$ for $n=1,2, \ldots$.
Each of the functions $e_n$ is locally constant;
in fact, if $\vert x-y \vert < n^{-1}$, then $e_n (x) =e_n (y)$.

\vskip 24pt

The following two simple results will turn out to be surprisingly useful,
and their proofs are nice illustrations of some of the ideas introduced in
this section.

\proclaim Lemma 3.11.  Suppose that $(E,\|\quad\|)$ is a normed vector space.  If $x_1, \ldots, x_n \in E$ then
$$\bigvee_{i=1}^n \|x_i\| = \|x_1\| \vee (\bigvee_{i=2}^n \|x_i - x_{i-1}\|)$$
and
$$\bigvee_{i=1}^n \|x_i\| = \|x_1 - \sum_{j=2}^n x_j\| \vee (\bigvee_{i=2}^n \|x_i\|).$$

\vskip 18pt
\noindent
{\it Proof.}  Consider the first equality.   It follows immediately from the ultrametric inequality that the right hand side is at most the left hand side.  Conversely, set $y_1 = x_1$ and $y_i = x_i - x_{i-1}$ for $2 \le i \le n$.  Then $x_i = \sum_{j=1}^i y_j$ for $1 \le i \le n$, and the reverse inequality also follows from the ultrametric inequality.  The proof of the second equality is similar and is left to the reader.
\hfil\break\rightline{\endproof}

\vskip 24pt

\proclaim Corollary 3.12. Let $e_i$, $1 \le i \le n$, denote the $i^{\rm th}$ coordinate vector in $K^n$.  Set $f_1 = e_1$ and $f_i = e_i - e_{i-1}$, $2 \le i \le n$.  Set $g_1 = e_1$ and $g_i = e_i - e_1$, $2 \le i \le n$.  Then both of the collections $f_1, \ldots, f_n$ and $g_1, \ldots, g_n$ are  orthornormal.

\vskip 18pt
\noindent
{\it Proof.}  We will prove the result for $f_1, \ldots, f_n$, and leave the proof for $g_1, \ldots, g_n$ to the reader.  Given $\alpha_1, \ldots, \alpha_n \in K$ we have, by Lemma 3.11, that
$$\eqalign{|\sum_i \alpha_i f_i| &= (\bigvee_{i=1}^{n-1} |\alpha_i - \alpha_{i+1}|) \vee |\alpha_n|\cr
&= \bigvee_{i=1}^n |\alpha_i|. \cr}$$  
\hfil\break\rightline{\endproof}

\vskip 24pt

   We will end this section with the analogue of the ``Cramer-Wold
device'' characterising the law of a random vector in terms of
the laws of all its ``one-dimensional projections''.
First we need the following fairly elementary piece of measure theory.

\noindent
\proclaim Lemma 3.13.
Suppose that $(E, \Vert~~~\Vert_E )$ is a separable Banach space with dual $E^*$.
Then $E^*$ generates the Borel $\sigma$-field of $E$.

\vskip 18pt
\noindent
{\it Proof.}.
It suffices to show that the map $x \mapsto \Vert x \Vert_E$ is measurable
with respect to the $\sigma$-field generated by $E^*$, but this follows by
standard arguments from the ultrametric analogue of the Hahn-Banach theorem
given in Appendix A.8 of [Schikhof, 1984] (cf. the proof of Proposition
7.1.1 in [Laha and Rohatgi, 1979]).
\hfil\break\rightline{\endproof}

\vskip 24pt
\noindent
\proclaim Lemma 3.14.
Let $(E, \Vert~~~\Vert_E )$ be a separable Banach space with dual $E^*$.
Let $(X_1 , \ldots , X_n )$ and $(Y_1 , \ldots , Y_n )$ be two $n$-tuples
of $E$-valued random variables.
Suppose that
$$
\bP[\chi( \sum_i T_i (X_i ))] =\bP[\chi ( \sum_i T_i (Y_i ))]
$$
for all $n$-tuples $(T_1 , \ldots , T_n )$ drawn from $E^*$.
Then
$(X_1 , \ldots , X_n )$ and $(Y_1 , \ldots , Y_n )$ have the same law.

\vskip 18pt
\noindent
{\it Proof.}.
Note that we can make $E^n$ into a Banach space with norm $\| \; \|_n$
by setting $\|(x_1, \ldots, x_n)\|_n = \bigvee_{i=1}^n \|x_i\|$.
The Borel $\sigma$-field of $E^n$ is then just $\cB^n$, where $\cB$ is
the Borel $\sigma$-field of $E$.  Each $(U_1, \ldots, U_n) \in (E^*)^n$
defines an element of $(E^n)^*$ via
$(x_1, \ldots, x_n) \mapsto \sum_i U_i(x_i)$,
and  all elements of $(E^n)^*$ arise this way for a unique choice
of $(U_1, \ldots, U_n)$.  The result of the lemma is thus just the
special case $n=1$ applied to the Banach space $E^n$, and so it suffices
to prove the $n=1$  case.

   Put $X=X_1$ and $Y=Y_1$.
Applying Lemma 3.13, we need to show  for any integer $m$ and any 
$S_1, \ldots,S_m \in E^*$ that the $m$-tuple
$(S_1(X), \ldots, S_m(X))$ has the same law as $(S_1(Y), \ldots, S_m(Y))$.
By Fourier uniqueness in $K^m$, it therefore
suffices to show that
$$
\bP[\chi( \sum_j \alpha_j S_j (X))] = \bP[\chi( \sum_j \alpha_j S_j (Y))]
$$
for all $(\alpha_1, \ldots, \alpha_m) \in K^m$;
but this follows from the hypothesis of the lemma with $T=\sum_j \alpha_j S_j$.
\hfil\break\rightline{\endproof}

\vfill\eject

\centerline{\bf 4. Properties of $K$-Gaussian probability measures}

\vskip 24pt

There are numerous (equivalent!) definitions of the class of Gaussian
 distributions
on $\bR$.   At one end of the spectrum, we have the explicit description
of Gaussian distributions in terms of their densities.   At the other
end we have characterisations such as the following, the form of
which goes back to Maxwell's derivation of the velocity
distribution for an ideal gas: 
a real-valued random variable $X$ has a centred Gaussian
distribution if and only if whenever
$X_1$, $X_2$ are two independent
copies of $X$ and $( \alpha_{11} , \alpha_{12} )$, $( \alpha_{21}, \alpha_{22} ) \in \bR^2$
are orthonormal, then $(X_1 , X_2 )$ has the same law as
$( \alpha_{11} X_1 + \alpha_{12} X_2 , \alpha_{21} X_1 + \alpha_{22} X_2 )$.

One can attempt to mimic such definitions when seeking an appropriate
definition for the class of Gaussian distributions on any real vector
space with suitable measurable structure.  The latter characterisation
extends in an obvious way to such general settings, whereas characterisations
such as the former require some further structure on the space (for example,
the existence of a Lebesgue measure) and they don't generalise as
easily or completely.
Work in the abstract theory of Gaussian measures (cf. 
[Fernique, 1975]) has shown that the latter definition is the
most ammenable to very general extension, and that such extensions typically
subsume extensions based on other characterisations.

As we have a notion of orthonormality in our local field setting,
the latter definition can also be mimicked there.  The
resulting theory is worked in the greatest possible generality
of ``measurable vector spaces over $K$'' in [Evans, 1989a].
Here we just content ourselves with the separable Banach space
case in order to streamline the presentation.

\vskip 24pt
\noindent
{\bf Definition 4.1.}
Let $E$ be a separable Banach space (over $K$) and suppose that $X$
is an $E$-valued random variable.
We say that the distribution of $X$ is a $K$-Gaussian probability measure (or, more simply, that $X$ is $K$-Gaussian) if when $X_1$, $X_2$ are two independent
copies of $X$ and $( \alpha_{11} , \alpha_{12} )$, $( \alpha_{21}, \alpha_{22} ) \in K^2$
are orthonormal, then $(X_1 , X_2 )$ has the same law as
$( \alpha_{11} X_1 + \alpha_{12} X_2 , \alpha_{21} X_1 + \alpha_{22} X_2 )$.

\vskip 24pt

Clearly, if $X=0$ almost surely, then $X$ is $K$-Gaussian.
Our first order of business obviously should be to show that there are
non-trivial $K$-Gaussian random variables.

\vskip 24pt
\noindent
\proclaim Theorem 4.2.
A $K$-valued random variable $X$ that is not almost surely 0 is $K$-Gaussian
if and only if the law of $X$ is normalised Haar measure restricted
to one of the $D$-modules $\rho^{-n} D$ for some $n \in \bZ$.  That is,
if and only if 
$$
\bP\{X \in dx\} =q^{-n} \Phi ( q^{-n} \vert x \vert )\lambda (dx)
= {1 \over {|\rho^{-n}|}} \b1_{\rho^{-n} D} (x) \lambda(dx)
$$
or, equivalently,
$$
\bP[\chi(\xi X)] = \Phi (q^n \vert \xi \vert ) .
$$

\vskip 18pt
\noindent
{\it Proof.}
Suppose that the distribution of $X$ has the given Fourier transform for some $n \in \bZ$.
From (2.4) and Fourier uniqueness, this is equivalent to the distribution of $X$ being as stated.
If $X_1 , X_2$ are two independent copies of $X$ and $( \alpha_{11} , \alpha_{12} )$,
$( \alpha_{21} , \alpha_{22} )$ is a pair of orthonormal vectors in $K^2$, then,
recalling (2.5),
$$
\eqalign{\bP[\chi(\xi_1&( \alpha_{11} X_1 + \alpha_{12} X_2 )+ \xi_2 ( \alpha_{21} X_1 + \alpha_{22} X_2 ))] \cr
& = \Phi (q^n \vert \xi_1 \alpha_{11} + \xi_2 \alpha_{21} \vert ) \Phi
(q^n \vert \xi_1 \alpha_{12} + \xi_2 \alpha_{22} \vert ) \cr
& = \Phi (q^n (\vert \xi_1 \alpha_{11} + \xi_2 \alpha_{21} \vert \vee \vert
\xi_1 \alpha_{12} + \xi_2 \alpha_{22} \vert) ) \cr
& = \Phi (q^n \vert \xi_1 ( \alpha_{11} , \alpha_{21} ) + \xi_2 ( \alpha_{21} , \alpha_{22} ) \vert ) \cr
& = \Phi (q^n (\vert \xi_1 \vert \vee \vert \xi_2 \vert) ) \cr
& = \Phi (q^n \vert \xi_1 \vert ) \Phi (q^n \vert \xi_2 \vert ) \cr
& = \bP[\chi(\xi_1 X_1) ]\bP[\chi(\xi_2 X_2)] . \cr}
$$

  From Fourier uniqueness, we have that $( \alpha_{11} X_1 + \alpha_{12} X_2 , \alpha_{21} X_1 + \alpha_{22} X_2 )$
has the same law as $(X_1 ,X_2 )$, and hence $X$ is $K$-Gaussian.

Conversely, suppose that $X$ is $K$-Gaussian.
Put $\varphi ( \xi ) =\bP[\chi(\xi X)]$.
Let $X_1$, $X_2$ be two independent copies of $X$.
Since $(1,1) , (0,1)$ is a pair of orthonormal vectors in $K^2$,
we have that $(X_1 +X_2 , X_2 )$ has the same law as $(X_1,X_2)$, and so
$( \varphi ( \xi ))^2 = \varphi ( \xi )$.
Thus $\varphi ( \xi ) \i \{ 0,1 \}$ for all $\xi \in K$.

Suppose that $\xi_0 \ne 0$ with $\varphi ( \xi_0 )=1$.
Such a $\xi_0$ must exist since $\varphi$ is continuous with $\varphi (0) =1$.
Note that $(1, \alpha) ,(0,1)$ is a pair of orthonormal vectors in $K^2$
 for each
$\vert \alpha \vert \le1$, and so $X_1 + \alpha X_2$ has the same law as $X_1$.
Thus $\varphi ( \xi ) \varphi ( \alpha \xi ) = \varphi(\xi)$ and so
$\varphi ( \alpha \xi_0 ) =1$.
This implies that $\varphi ( \xi ) =1$ for all $\vert \xi \vert \le \vert \xi_0 \vert$.
Now $\varphi \ne 1$, since $X$ is not almost surely $0$; so we must have that $\varphi ( \xi ) = \Phi (q^n \vert \xi \vert )$ for some $n \in \bZ$.
\hfil\break\rightline{\endproof}

\vskip 24pt
\noindent
\proclaim Corollary 4.3.
If $X$ is a $K$-valued $K$-Gaussian random variable, then $X \in \cL^{\infty} (\bP)$.
We have that
$$
\bP[\chi(\xi X)] = \Phi ( \Vert X \Vert_{\infty} \vert \xi \vert ) .
$$
If $X$ is not almost surely 0, then
$$
\bP\{X \in dx \} = \Vert X \Vert_{\infty}^{-1} \Phi ( \Vert X \Vert_{\infty}^{-1}
\vert x \vert ) \lambda (dx) .
$$

\vskip 18pt
\noindent
{\it Proof.}
Clear.
\hfil\break\rightline{\endproof}

\vskip 24pt 
\noindent
{\bf Remark.} Our definition of  $K$-Gaussian probability measures
 is the analogue
of a definition for centred real Gaussian probability measures.
We could, of course, also consider the ``non-centred'' case by
introducing random variables of the form $x+X$, where $X$ is
$K$-Gaussian and $x \in E$.  Note from Corollary 4.3 that if
$X$ is $K$-valued and $|y-z| \le \|X\|_\infty$, then $y+X$ has
the same distribution as $z+X$, and so the ``shift parameter''
$x$ is not uniquely determinable from the distribution of $x+X$.

\vskip 24pt

\noindent
\proclaim Theorem 4.4.
Suppose that $E$ is a separable Banach space.
If $X$ is an $E$-valued random variable, then $X$ is $K$-Gaussian if and
only if $T(X)$ is $K$-Gaussian for all $T \in E^*$.

\vskip 18pt
\noindent
{\it Proof.}
Suppose first of all that $T(X)$ is $K$-Gaussian for all $T \in E^*$.
Let $X_1$, $X_2$ be two independent copies of $X$.
Fix an orthonormal pair of vectors $( \alpha_{11}, \alpha_{12})$,
$( \alpha_{21} , \alpha_{22})$ and a pair of functionals $T_1, T_2 \in E^*$.
From Lemma 3.14, 
in order to show that $X$ is $K$-Gaussian, we need to check that
$$
\leqalignno{\bP[\chi(T_1&( \alpha_{11} X_1 + \alpha_{12} X_2 )+T_2 ( \alpha_{21} X_1 + \alpha_{22} X_2 ))] \cr
&=\bP[\chi(T_1 (X_1 )+T_2 (X_2 ))] .&(4.4.1) \cr}
$$

  From Corollary 4.3 we see that the left-hand side of (4.4.1) is given by
$$
\Phi ( \Vert ( \alpha_{11} T_1 + \alpha_{21} T_2 ) X \Vert_{\infty} \vee
\Vert ( \alpha_{12} T_1 + \alpha_{22} T_2 ) X \Vert_{\infty} ) .
$$
Suppose that the dimension of the span of $\{ T_1 (X) ,T_2 (X) \}$
in $\cL^\infty(\bP)$ is  $2$ (the cases where the
dimension is 1 or 0 can be handled similarly and more easily).
Let $Y_1$, $Y_2$ be an orthonormal basis for the span of
 $\{ T_1 (X) ,T_2 (X) \}$
in $\cL^{\infty} (\bP)$  and write
$$
T_1 (X) = \beta_{11} Y_1 + \beta_{12} Y_2 , 
$$
$$
T_2 (X) = \beta_{21} Y_1 + \beta_{22} Y_2 . 
$$
We have
$$
\eqalign{\Vert ( \alpha_{11} T_1 + \alpha_{21} T_2 ) X \Vert_{\infty}&= \Vert
\alpha_{11} ( \beta_{11} Y_1 + \beta_{12} Y_2 )+ \alpha_{21} ( \beta_{21} Y_1 + \beta_{22} Y_2 ) \Vert_{\infty} \cr
&= \vert \alpha_{11} \beta_{11} + \alpha_{21} \beta_{21} \vert \vee \vert
\alpha_{11} \beta_{12} + \alpha_{21} \beta_{22} \vert . \cr}
$$
Similarly,
$$
\Vert ( \alpha_{12} T_1 + \alpha_{22} T_2 ) X \Vert_{\infty} = \vert
\alpha_{12} \beta_{11} + \alpha_{22} \beta_{21} \vert \vee \vert
\alpha_{12} \beta_{12} + \alpha_{22} \beta_{22} \vert .
$$
One can now readily check that the left-hand side of (4.4.1) is just
$$
\eqalign{\Phi ( \vert \beta_{11}&( \alpha_{11} , \alpha_{12} ) + \beta_{21}
( \alpha_{21} , \alpha_{22} ) \vert \cr
&\vee \vert \beta_{12} ( \alpha_{11} , \alpha_{12} ) + \beta_{22} ( \alpha_{21} , \alpha_{22} ) \vert ) \cr
&= \Phi (( \vert \beta_{11} \vert \vee \vert \beta_{21} \vert ) \vee
( \vert \beta_{12} \vert \vee \vert \beta_{22} \vert )) \cr}
$$
by the orthonormality of $( \alpha_{11} , \alpha_{12} )$, $( \alpha_{21} , \alpha_{22} )$.

Also, the right-hand side of (4.4.1) is just
$$
\eqalign{\Phi ( \Vert \beta_{11} Y_1&+ \beta_{12} Y_2 \Vert_{\infty} \vee
\Vert \beta_{21} Y_1 + \beta_{22} Y_2 \Vert_{\infty} ) \cr
&= \Phi (( \vert \beta_{11} \vert \vee \vert \beta_{12} \vert ) \vee (
\vert \beta_{21} \vert \vee \vert \beta_{22} \vert )) . \cr}
$$
Therefore (4.4.1) holds and $X$ is $K$-Gaussian.

   Suppose now for the converse that $X$ is $K$-Gaussian and $T \in E^*$. 
As above,
let $X_1,X_2$ be two independent copies of $X$, and let 
 $( \alpha_{11}, \alpha_{12})$,
$( \alpha_{21} , \alpha_{22})$
be a pair of orthonormal vectors.
By definition, the pair
$(( \alpha_{11} X_1, \alpha_{12} X_2)$, $( \alpha_{21} X_1, \alpha_{22} X_2))$
has the same law as $(X_1, X_2)$; and so, by the linearity of $T$,
the pair
$(( \alpha_{11} T(X_1), \alpha_{12} T(X_2))$, $( \alpha_{21} T(X_1), \alpha_{22} T(X_2)))$ has the same law as $(T(X_1), T(X_2))$.  As $T(X_1), T(X_2)$ are
two independent copies of $T(X)$, it follows from the definition that $T(X)$
is $K$-Gaussian.
\hfil\break\rightline{\endproof}

\vskip 24pt
\noindent
\proclaim Corollary 4.5.
Suppose that $E$ is a separable Banach space.
If $(X_n )_{n=0}^{\infty}$ is a sequence of $E$-valued $K$-Gaussian
random variables such that $X$ converges in distribution
 as $n \rightarrow \infty$
to some random variable $X$, then $X$ is also $K$-Gaussian.

\vskip 18pt
\noindent
{\it Proof.}
By Theorem 4.4, we need only check that $T(X)$ is $K$-Gaussian for all
$T \in E^*$.
However,
$$
\eqalign{\bP[\chi(\xi T(X))]&=\bP[\chi( \lim_{n \rightarrow \infty} \xi T(X_n ))] \cr
&= \lim_{n \rightarrow \infty} \bP[\chi( \xi T(X_n ))] \cr
&= \lim_{n \rightarrow \infty} \Phi ( \Vert T(X_n ) \Vert_{\infty} \vert \xi \vert ) \cr
&= \Phi ( \Vert T(X) \Vert_{\infty} \vert \xi \vert ) ,\cr}
$$
and the result follow from Corollary 4.3.
\hfil\break\rightline{\endproof}

\vskip 24pt

   In Theorem 4.2 we saw that the set of
 possible laws of $K$-valued, $K$-Gaussian
random variables is just the set of normalised Haar measures on
compact $D$-submodules of $E$.  The following two results, Theorem 4.6
and Theorem 4.7, show that the verbatim extension of this result to
Banach space valued random variables holds.

\vskip 24pt
\proclaim Theorem 4.6. Suppose that $E$ is a separable Banach space
and let $X$ be a $K$-Gaussian random variable.
Set
$$
S=\{x \in E : |T(x)| \le \Vert T(X) \Vert_\infty,\quad \forall T \in E^*\}.
$$
i) The set $S$ is a $D$-submodule.\hfil\break
ii) If $x \in S$, then the law of $x+X$ coincides with the law of $X$.
Otherwise, the law of $x+X$ and the law of $X$ are mutually singular.
\hfil\break
iii) The set $S$ is the closed support of the law of $X$.
\hfil\break
iv) The set $S$ is compact.\hfil\break
v) The law of $X$ is normalised Haar measure on $S$.\hfil\break
vi) If $M$ is a measurable vector subspace of $E$ then $\bP\{X \in M\}$ is either 1 or 0,
depending on whether or not $S \subset M$.

\vskip 18pt
\noindent
{\it Proof.} i) This is immediate from properties (3.2)-(3.4).

\noindent
ii) From Theorem 4.4, Theorem 4.2 and Corollary 4.3, the law of $T(X)$
for $T \in E^*$ is normalised Haar measure restricted to the $D$-module
$\{y \in K : |y| \le \|T(X)\|_\infty$.
Hence, if $x \in S$ we have that the law of $T(x+X) = T(x)+T(X)$
coincides with the law of $T(X)$  for all $T \in E^*$.
Thus, by Lemma 3.14, the law of $x+X$ coincides with that of $X$.

   If $x \notin S$ then there exists $T \in E^*$ such that $|T(x)| > \Vert T(X) \Vert_\infty$.
Then, by the isosceles triangle property, $|T(x+X)|=|T(x)|$,
almost surely, and hence the law of $x+X$ and the law of $X$ are mutually 
singular.

\noindent
iii)  It is clear that $S$, as an intersection of closed sets, is closed.
Let $\{T_i\}_{i=1}^\infty$ be a countable dense subset of $E^*$.
Then
$$
\bP\{X \in S\}=\bP(\bigcap_{i=1}^\infty \{ |T_i(X)| \le \Vert
T_i(X)\Vert_\infty\}) = 1.
$$
Conversely, suppose that $x \in S$ and $U$ is an open neighbourhood of $x$.
Let $\{x_i\}_{i=1}^\infty$ be a countable dense subset of $S$.
Then $\bigcup_{i=1}^\infty  [x_i + (U-x)]$ covers $S$ and hence
$\bP\{ X \in x_i + (U -x)\} > 0$ for at least one $i$;
but, by (ii), $\bP\{X \in U\} = \bP\{X \in (x_i - x) + U\}$, since
$x_i - x \in S$.

\noindent
iv)  As $E$ is complete and separable, all probability measures on $E$
are tight and so there exists a compact set $C \subset S$ such that 
$\bP\{X \in C\}>0$.
Let $G$ be the smallest closed $D$-submodule  containing $C$.
We claim that $G$ is also compact.
Given $\epsilon>0$, there exists a finite set $\{x_1^\epsilon,\ldots,x_{n(\epsilon )}^\epsilon\}\subset C$
such that if $x \in C$ then $\Vert x -x_i^\epsilon\Vert < \epsilon$
for some $x_i^\epsilon$.
The smallest closed $D$-submodule containing $\{x_1^\epsilon,\ldots,x_{n(\epsilon )}^\epsilon\}$
 is $G^\epsilon = (D x_1^\epsilon )+\ldots+(D x_{n(\epsilon )}^\epsilon )$.
Clearly, $G^\epsilon$ is compact.
Moreover, from the ultrametric inequality it is clear that if $x \in G$,
then there exists $y \in G^\epsilon$ such that $\Vert x-y \Vert < \epsilon$.
Thus $G$ is totally bounded and hence compact.

   Part (iv) will now follow if $G$ has only finitely many distinct cosets in $S$;
but this must be the case, since otherwise we could find infinitely many
disjoint cosets $G_1, G_2 ,...$ for which, by part (ii),
$\bP\{X \in G_i \} = \bP\{X \in G\} > 0$.

\noindent
v) This is immediate from parts (i)-(iv).

\noindent
vi)  We begin by showing that $\bP\{X \in M\}$ is either $0$ or $1$.
Let $X_1 , X_2$ be two independent copies of $X$.
For $n=1,2, \ldots$, set
$$
A_n = \{ X_1 + \rho^n X_2 \in M , \rho^n X_1 + X_2 \not\in M \} .
$$

   We claim first of all that the vectors $(1, \rho^n )$, $( \rho^n ,1)$ are
an orthonormal pair.
Clearly, $\vert (1, \rho^n ) \vert = \vert ( \rho^n ,1) \vert =1$.
From the ultrametric inequality and symmetry, we need only show that there
is no $\alpha \in K$ such that
$$
\vert (1, \rho^n ) + \alpha ( \rho^n ,1) \vert = \vert 1+ \alpha \rho^n \vert
\vee \vert \rho^n + \alpha \vert <1 .
$$
If this was so, then $\vert 1+ \alpha \rho^n \vert <1$, and we see from
the isosceles triangle property that $\vert \alpha \rho^n \vert =1$
and hence $\vert \alpha \vert =q^n$.
Again applying the isosceles triangle property, this implies that
$\vert \rho^n + \alpha \vert =q^n >1$, which is a contradiction.

By definition we therefore have that $(X_1 + \rho^n X_2 , \rho^n X_1 +X_2 )$
has the same law as $(X_1 , X_2 )$, and hence
$$
\bP(A_n ) =\bP\{X \in M\}\bP\{X \not\in M\}.
$$

Observe also that if $m \ne n$, then the matrix
$$
\left(\matrix{1&\rho^m\cr
1&\rho^n\cr}\right)
$$
is invertible.
Thus, if both $X_1 + \rho^m X_2$ and $X_1 + \rho^n X_2$
belong to $M$, then both $X_1$ and $X_2$ (and hence both $\rho^m X_1 + X_2$ and $\rho^m X_1 + X_2$) belong to $M$, and thus
$A_m \cap A_n = \emptyset$.
It follows that
$$
1 \ge \bP(\bigcup_n A_n ) =\sum_n \bP(A_n )= \sum_n\bP\{X \in M\}\bP\{X \not\in M\},
$$
and so  we conclude that $\bP\{X \in M\} \in \{ 0,1 \}$.

   Now, if $S \subset M$ then it follows from part (iii) that $\bP\{X \in M\} = 1$.
Conversely, suppose that $\bP\{X \in M\} = 1$.
If there exists $x \in S$ such that $x \notin M$ then $M$ and $x+M$ are disjoint;
but this is impossible, since $\bP\{X \in x+M\} = \bP\{X \in M\} = 1$ by 
part (ii).
\hfil\break\rightline{\endproof}

\vskip 24pt

   The $D$-submodule $S$ is, in many ways, analogous to the reproducing
kernel Hilbert space (rkhs) in the theory of (ordinary) Gaussian 
random variables.  There,
``Cameron-Martin'' type theorems (cf. [Feldman, 1958]
or [Hajek, 1959]) state that a shift of a Gaussian random variable
by an element of the rkhs gives a random variable with an equivalent distribution, whereas a shift by
an element not in the rkhs gives a random variable with a mutually singular
distribution.  Moreover, the support of the law of the Gaussian
random variable is the closure of the rkhs, and the rkhs uniquely determines
the law of the Gaussian random variable. 

 However, although it
is known that the probability a Gaussian random variable belongs to
a measurable vector subspace is either $0$ or $1$ (see Theorem 1.2.1
of [Fernique, 1975]), there does not seem to be a counterpart to
part (vi) of Theorem 4.6 giving a condition that determines
which branch of the dichotomy holds. Also, in
certain Gaussian situations
 is possible to  obtain a similar zero-one law for the probability
of belonging to a measurable
subgroup rather than a vector subspace (cf. [Kallianpur, 1970],
[Jain, 1971] and [Cambanis and Rajput, 1973]),
and it is clear that such a result will not hold in general for our setting.
For if $X$ is a $K$-valued, $K$-Gaussian random variable with $\Vert X \Vert_{\infty} >0$,
then $G= \{ x: \vert x \vert \le q^{-1} \Vert X \Vert_{\infty} \}$ is a
subgroup of $K$ and $0<\bP\{X\in G\}=q^{-1} <1$.

   The converse to Theorem 4.6 holds.

\vskip 24pt

\proclaim Theorem 4.7.
Suppose that $E$ is a separable Banach space and that
$G$ is a compact $D$-submodule of $E$.
Suppose that $X$ is an $E$-valued random variable
with law that is normalised Haar measure on $G$.
Then $X$ is a $K$-Gaussian random variable for which, in the notation
of Theorem 4.6, $S=G$.

\vskip 18pt
\noindent
{\it Proof.}
In order to show that $X$ is $K$-Gaussian, it suffices by Theorem 4.4
to show that $T(X)$ is $K$-Gaussian for all $T \in E^*$.
As $T$ is continuous, it follows both that the support of the law of $T(X)$
is the set $T(G)$ and that this set is compact.  Further, as
$T$ is linear, it follows both that $T(G)$ is a $D$-submodule of $K$,
and that the law of $y + T(X)$ coincides with the law of $T(X)$ whenever
$y \in T(G)$.  Thus the law of $T(X)$ is normalised Haar measure on
the compact $D$-submodule $T(G)$.
Since the only compact $D$-submodules of $K$ are sets of the 
form $\rho^n D$ for some
$n \in \bZ$ and $\{0\}$,  Theorem 4.2 gives that $T(X)$
is $K$-Gaussian.
Part (iii) of Theorem 4.6 shows that $S=G$.
\hfil\break\rightline{\endproof}

\vskip 24pt
  We finish this section with some results which further 
reinforce the connection
between the properties of $K$-Gaussian random variables and the local
field theory of orthogonality.  One might imagine that as
Theorem 4.6 and Theorem 4.7 appear to be
the definitive ``structure'' results, all our proofs from now on
will appeal to them.  However, as the proofs of the following results
show, it is often easier to appeal to Theorem 4.4 and the one
dimensional characterisations of Theorem 4.2 and Corollary 4.3.

\vskip 24pt
\noindent
\proclaim Theorem 4.8.
i) If $X=(X_1 ,\ldots ,X_n )$ is a $K^n$-valued $K$-Gaussian random variable,
then $\{X_1 , \ldots , X_n \}$ is an orthogonal set in $\cL^{\infty} (\bP)$ if
and only if $X_1 , \ldots , X_n$ are independent.\hfil\break
ii) If $X_1, \ldots, X_n$ are independent $K$-valued, $K$-Gaussian
random variables, then the random vector\hfil\break 
$(X_1, \ldots, X_n)$ is
$K$-Gaussian.

\vskip 18pt
\noindent
{\it Proof.}
i) Suppose that $\{X_1 , \ldots , X_n \}$ is an orthogonal set.
From Theorem 4.4 we have that $\xi \cdot X$ is $K$-Gaussian for all
$\xi = ( \xi_1 , \ldots , \xi_n ) \in K^n$.
Applying Corollary 4.3, we have
$$
\eqalign{\bP[\chi( \xi \cdot X)]&=\bP[\chi(1( \xi \cdot X))] \cr
&= \Phi ( \Vert \xi \cdot X \Vert_{\infty} ) \cr
&= \Phi ( \Vert \sum_i \xi_i X_i \Vert_{\infty} ) \cr
&= \Phi (\bigvee_{i=1}^n \vert \xi_i \vert \Vert X_i \Vert_{\infty} ) \cr
&= \prod_{i=1}^n \Phi ( \Vert X_i \Vert_{\infty} \vert \xi_i \vert ) \cr
&= \prod_{i=1}^n \bP[\chi( \xi_i X_i )] , \cr}
$$
and the result follows from Fourier uniqueness.

The proof of the converse essentially consists of reversing the
above chain of inequalities, and we omit it.

\noindent
ii) From Corollary 4.3 we know that the law of $(X_1, \ldots, X_n)$ is just normalised Haar
measure on the set $\prod_{i=1}^n \{x \in K : |x| \le \|X_i\|_\infty\}$.
As this set is certainly a compact $D$-submodule of $K^n$, we conclude
from Theorem 4.7 that $(X_1, \ldots, X_n)$ is $K$-Gaussian. 
\hfil\break\rightline{\endproof}

\vskip 24pt
\proclaim Theorem 4.9.
Suppose that $E$ and $F$ are separable Banach spaces, $A:E \rightarrow F$
is a continuous linear operator, and $X$ is an $E$-valued, $K$-Gaussian
random variable.  Then $A(X)$ is an $F$-valued, $K$-Gaussian
random variable.

\vskip 18pt
\noindent
{\it Proof.} This is immediate from Theorem 4.4 once we observe that if
$T \in F^*$, then $T \circ A \in E^*$.
\hfil\break\rightline{\endproof}

\vskip 24pt
\noindent
\proclaim Corollary 4.10.
Suppose that $X=(X_1 , \ldots , X_n )$ is a $K^n$-valued random variable.
Then $X$ is $K$-Gaussian if and only if for some $m$ there exists a vector
$Y=(Y_1 , \ldots , Y_m )$ of independent $K$-valued $K$-Gaussian random
variables, and an $m \times n$ matrix $A$ such that $X=YA$.

\vskip 18pt
\noindent
{\it Proof.}
Suppose that $X=YA$ with $Y$ and $A$ as above.
From part (ii) of Theorem 4.8 we have that $Y$ is $K$-Gaussian,
and Theorem 4.9 gives that $X$ is also $K$-Gaussian.

Conversely, suppose that $X$ is $K$-Gaussian.
Let $\{ Y_1 , \ldots , Y_m \}$ be an orthonormal basis for the linear
span of
$\{ X_1 , \ldots , X_n \}$ in $\cL^{\infty} (\bP)$.
Since $Y=(Y_1 , \ldots , Y_m ) =XB$ for some $n \times m$ matrix $B$,
we have by Theorem 4.9 that $Y$ is $K$-Gaussian.
From part (i) of Theorem 4.8, $Y_1 , \ldots ,Y_m$ are independent and the result follows.
\hfil\break\rightline{\endproof}

\vskip 24pt
\noindent
\proclaim Theorem 4.11.
Suppose that $E$ satisfies is a separable Banach space
 and  $X=(X_1 , \ldots , X_n )$ is a vector of independent, identically
distributed, $E$-valued, $K$-Gaussian random variables.
Suppose that $\{ \alpha_1 , \ldots , \alpha_m \}$ is an orthogonal set in $K^n$.
Define an $n \times m$ matrix by $A = ( \alpha_1' , \ldots , \alpha_m' )$.
Then $Y=XA$ is a vector of independent $K$-Gaussian random variables.

\vskip 18pt
\noindent
{\it Proof.}  Applying Theorem 4.4 we reduce to the case $E=K$.
For $\xi =( \xi_1 , \ldots , \xi_m ) \in K^m$ we have, setting
$\sigma = \Vert X_1 \Vert_{\infty}$, that
$$
\eqalign{\bP[\chi( \xi \cdot Y)]&=\bP[\chi(XA \xi' )] \cr
&= \bP[\chi(X \cdot ( \xi A' ))] \cr
&= \prod_{j=1}^n \Phi ( \sigma \vert ( \xi A' )_j \vert ) \cr
&= \Phi ( \sigma \bigvee_{j=1}^n \vert ( \xi A' )_j \vert ) \cr
&= \Phi ( \sigma \vert ( \xi A' ) \vert ) \cr
&= \Phi ( \sigma \vert \sum_{i=1}^m \xi_i \alpha_i \vert ) \cr
&= \Phi ( \sigma \bigvee_{i=1}^m \vert \xi_i \vert \vert \alpha_i \vert ) \cr
&= \prod_{i=1}^m \Phi ( \sigma \vert \alpha_i \vert  \vert \xi_i \vert ) , \cr}
$$
and the result follows from Fourier uniqueness and Theorem 4.2.
\hfil\break\rightline{\endproof}

\vskip 24pt
\noindent
{\bf Remark.}
One of the main reasons why the (real) Gaussian theory is so
tractable is that  law of a centred Gaussian random vector
 is completely determined by its
family of 2-dimensional marginal distributions or, more precisely, by its
covariance matrix.
There is, however, no fixed integer $n$ such that the law of every $K$-Gaussian
random vector is determined by its family of $n$-dimensional marginal distributions.
Suppose that $\{ Z_1 , \ldots , Z_{n+1} \}$ is a set of independent
identically distributed $K$-valued, $K$-Gaussian random variables with
$\Vert Z_1 \Vert_{\infty} =1$.  The random vector $(Z_1, \ldots Z_{n+1})$
is $K$-Gaussian by part (ii) of Theorem 4.8.
Set $X=(Z_1 , \ldots , Z_n , Z_{n+1} )$ and
$Y=(Z_1 , \ldots , Z_n , Z_1 + \ldots + Z_n )$.
One can check, using Theorem 4.11, that for each set of indices
$\{ i_1 , \ldots , i_n \} \subset \{ 1, \ldots ,n+1 \}$ we have that
$(X_{i_1} , \ldots , X_{i_n} )$ has the same law as $(Y_{i_1} , \ldots , Y_{i_n} )$.
Clearly, however, $X$ does not have the same law as $Y$.

\vskip 24pt
   The following calculations will be used in 
\S 5 and are   immediate from Theorem
4.4, part (ii) of Theorem 4.8, and Corollary 3.12.

\vskip 24pt

\proclaim Corollary 4.12.  Suppose that  $E$ is a separable Banach space
and $X_1, \ldots, X_n$ are independent, identically distributed, $E$-valued $K$-Gaussian random variables.  Then the random variables $X_1, X_2 - X_1, \ldots, X_n - X_{n-1}$ are independent and identically distributed random variables with the same common distribution as  $X_1, \ldots, X_n$.  The same is also true for the collection $X_1, X_2 - X_1, \ldots, X_n - X_1$.

\vfill\eject

\vskip 30pt

\centerline{\bf 5. Construction and first properties of K-Brownian motion}

\vskip 24pt

   One of the main reasons for the central importance
of (real) Brownian motion in the theory of stochastic processes is that
it can be defined in numerous ways that at first sight seem to bear
little relation to each other.  For instance, Brownian motion on the line
is characterised as:

\vskip 12pt\noindent
i) the centred Gaussian process with continuous paths and covariance kernel 
$(s,t) \mapsto s \wedge t$;\hfil\break
ii) the centred Gaussian process with continuous paths and 
reproducing kernel Hilbert
space the Cameron-Martin space of absolutely continuous functions 
$f:[0,1] \rightarrow \bR$ for which $f(0)=0$ and
$\int_0^1 |f'(t)|^2 \, dt < \infty$ 
(here we are thinking of Brownian motion as a process indexed by
$[0,1]$ rather than the whole half-line);\hfil\break
iii) the process with continuous paths and centred, stationary, independent increments;\hfil\break
iv) the Markov process with continuous paths and generator extending 
${1 \over 2}{{d^2} \over {dx^2}}$;\hfil\break
v) the continuous martingale with quadratic variation process $t \mapsto t$.

\vskip 12pt\noindent
Brownian motion with values in $\bR^d$ has similar characterisations.
Analogues of a subset of conditions (i) - (v) hold for the
various possible multiparameter relatives of Brownian motion such as
the Brownian sheet and L\'evy's multiparameter Brownian motion.

   Our aim in this section is to use Brownian motion
and its multiparameter relatives as guides for where to look
for interesting and canonical processes in our local field setting.  Because
we have the $K$-Gaussian theory at our disposal and we
know what object plays the role of the reproducing kernel
Hilbert space in that theory (see the remarks after Theorem 4.6), 
our definition will
be a counterpart of characterisation (ii) above (informed by the
wisdom that interesting local field objects often come from taking
a corresponding real object and ``replacing $\cL^2$ by $\cL^\infty$''.)

\vskip 24pt
\noindent
{\bf Definition 5.1.} Let $N$ and $d$ be positive integers.  A
$(N,d)$ $K$-Brownian motion is a $C(D^N, K^d)$-valued, $K$-Gaussian random variable  $X$ such that the closed support in $C(D^N,K^d)$ of the distribution of $X$ is the compact $D$-submodule $S$ consisting of functions $f$ such that $|f(t)| \le 1$ for all $t \in D^N$ and $|f(s) - f(t)| / |s-t| \le q^{-1}$ for all $s \ne t \in D^N$.

\vskip 24pt
\noindent
{\bf Remarks.} i) It is clear that the set $S$ is a closed $D$-module,
 and, as $S$ is equicontinuous, it follows from the Arzela-Ascoli 
theorem that $S$ is compact.\hfil\break
ii) The index set of $X$ is $D^N$, rather than $[0,1]^N$
or $\bR_+^N$.  As $K^d$ is totally disconnected, the only continuous maps
from $\bR_+^N$ to $K^d$ are the constants.  There are certainly interesting
discontinuous $K^d$-valued processes indexed by $\bR_+$ (see [Evans, 1989b]
or [Albeverio and Karwowski, 1991, 1994]).\hfil\break
iii) The reader might find it troubling that the definition involves
difference ratios rather than some sort of derivative.  
An explanation of why this is appropriate would lead us too far
afield, and we refer the reader to \S\S 26 and 27 of [Schikhof, 1984]
for an indication of why this approach is natural.\hfil\break
iv) Whilst we this definition has the right ``feel'', we don't
expect the reader immediately to find it particularly compelling.  
Rather, ``the
proof of the pudding will be in the eating'', as we show that the
object we propose exhibits many of the other features of Brownian motion.
For example, we will obtain analogues of the Markov and strong Markov properties
(Lemmas 5.2 and 5.3, and Corollary 6.2)
 and the independent increments property (Lemma 5.4).  
Further, we will show   in \S\S 6 and 7 that
we have a theory of ``additive functionals'' and a ``probabilistic
potential theory'' much like that of Brownian motion.
Of course, the differences between the structure of the real numbers and
a local field are so great that we shouldn't expect there to be a perfect
match between the two theories.  For example, the properties (i) and (v)
don't even seem to have natural translations into the local field world.

\vskip 24pt

  Although Definition 5.1 is particularly succint, it is not
easy to work with when trying to deduce or prove results about a
$K$-Brownian motion $X$.  For this purpose, the following ``bare hands''
construction will be much more useful.

   Before giving this construction of $X$ in full detail, we describe it
somewhat informally.  First consider the case $d=1$.  
Recall from \S 2
that we can think of the sub-balls of $D$ as being arrayed
in an infinite, rooted, $q$-ary tree.
A similar picture holds for the sub-balls of $D^N$, except
that now we have a $q^N$-ary tree.  
The root is $D^N$ itself, 
the nodes at level $k$ are the balls of radius
$q^{-k}$ (= cosets of $\rho^k D^N$), and the $q^N$ ``children'' 
of such a ball
are the $q^N$ balls of radius $q^{-(k+1)}$ that it contains.
We  attach a $K$-valued, $K$-Gaussian ``weight'' to each
such ball.  The weights are independent, and the weight
assigned to a ball has $\cL^\infty(\bP)$ norm equal to the  radius
of the ball.  Each point in $t \in D^N$ is associated with the unique
sequence of balls of decreasing radius that contain it, and we obtain
$X(t)$ by summing the attached weights.  

   Thus $X(t)$ is a sum of
weights that are orthogonal in $\cL^\infty(\bP)$ (see Theorem 4.8) and have norms 
$1, q^{-1}, q^{-2}, \ldots$, and hence $\|X(t)\|_\infty = 1$.  Similarly,
if $|s-t| = q^{-k}$, then $s$ and $t$ are contained in the same ball of
radius $q^{-k}$ but in different balls of radius $q^{-(k+1)}$.
Thus $X(s) - X(t)$ is the sum of the weights attached to balls containing
$s$ that are of radius at most $q^{-(k+1)}$ minus the 
sum of the weights attached to balls containing
$t$ that are of radius at most $q^{-(k+1)}$.  Hence $X(s) - X(t)$ is a sum
of random variables that are orthogonal in $\cL^\infty(\bP)$ and have norms
 $q^{-(k+1)}, q^{-(k+1)}, q^{-(k+2)}, q^{-(k+2)}, \ldots$.  Consequently,
$\|X(s) - X(t)\|_\infty = q^{-(k+1)} = q^{-1} |s-t|$.  These
two calculations indicate that we are on the right track to building
an object that satisfies Definition 5.1.

   The construction for general $d$ is effected by identifying
$C(D^N, K^d)$ with $(C(D^N,K))^d$ and building each $C(D^N,K)$-valued
component as an independent copy of the random variable constructed
in the $d=1$ case.  In  more concrete terms, this construction is much the same
as the one above, except that the weight attached to a ball of
radius $q^{-k}$ is now a vector of $d$ i.i.d. $K$-valued, $K$-Gaussian
random variables, each with $\cL^\infty(\bP)$ norm $q^{-k}$.

   We will now formalise this construction.
Let $\cD(N;n)$ denote the balls of radius $q^{-n}$ that are contained
in $D^N$, and set $\Gamma(N) = \cup_{n=0}^\infty \cD(N;n)$.       
Put  $\Omega=\prod_{C\in\Gamma(N)}\prod_{i=1}^d K_{C,i}$, where $K_{C,i}$ is a copy of $K$.   Let $Z_{C,i}$ denote the coordinate projection from $\Omega$ onto $K$ and write $Z_C$ for the map from $\Omega$ onto $K^d$ given by $Z_C=(Z_{C,i})_{i=1}^d$.   The $K^d$-valued random variable $Z_C$
will be the ``weight'' attached to the ball $C$.
Let $\zeta_{C,i}$ denote the $K$-Gaussian probability measure on $K_{C,i}$ for which the identity map has $\cL^\infty$ norm $q^{-n}$.  That is, $\zeta_{C,i}$ is the Haar measure $\lambda $ restricted to the copy of $\rho^n D$ in $K_{C,i}$ and renormalised to have total mass $1$.  Write $\zeta=\prod_{C\in\Gamma(N)}\prod_{i=1}^d \zeta_{C,i}$,
so that the random variables $Z_{C,i}$, $C \in \Gamma(N)$, $1 \le i \le d$, are
independent under $\zeta$.
  Let $\bP$ be the completion of $\zeta$ and let $\cG$ be the corresponding completion of the product $\sigma$-field on $\Omega$.

   Put
$$\Omega^*=\bigcap_{n=0}^\infty \bigcap_{C \in \cD(N;n)} \bigcap_{i=1}^d \{\omega:|Z_{C,i}(\omega)| \le q^{-n} = \|Z_{C,i}\|_\infty\},$$
so that $\bP(\Omega^*)=1$.  For a positive integer $n$ set
$$X^n(\omega)=\cases{\sum_{k=0}^n \sum_{C \in \cD(N;k)} Z_C(\omega) \b1_C, &if $\omega\in \Omega^*$, \cr
                     0,                                  &if $\omega\notin \Omega^*$;\cr}$$
so that, essentially, we build the value of $X^n$ at the index $t \in D^N$
by only adding up the weights attached to balls contain $t$ that have
radius $q^{-n}$ or larger.

  It follows from part (ii) of Theorem 4.8 and Theorem 4.9 that $X^n$
is a $C(D^N, K^d)$-valued, $K$-Gaussian random variable.  By definition
of $\Omega^*$ and the ultrametric inequality, we have 
$\|X^n(\omega)\|_C \le 1$  and 
$\|X^m(\omega) - X^n(\omega)\|_C \le q^{-(m+1)}$, $m < n$ 
for all $\omega \in \Omega$.  Thus $\{X^n\}_{n=1}^\infty$ converges
to a $C(D^N, K^d)$-valued random variable that we will denote by $X$.
By Corollary 4.5, $X$ is $K$-Gaussian.

\vskip 24pt
\noindent
{\bf Note:} {\it Until further notice, the notation $X$ will be reserved specifically for the  $C(D^N, K^d)$-valued, $K$-Gaussian random variable constructed as above.}

\vskip 24pt

   We need to show that $X$ is, in fact, a $(N,d)$
$K$-Brownian motion, and we will do so at the end of this section.  First,
however, we need to develop some elementary properties of $X$.

   Specify a family $\{\theta_t:t \in D^N\}$ of measure preserving bijections on $\Omega$ by requiring that $Z_{C,i}(\theta_t(\omega))=Z_{t+C,i}(\omega)$.  Observe that $\theta_0$ is the identity map and $\theta_s \circ \theta_t = \theta_{s+t}$, so that this family forms a group under the composition operation.  Note that $X(\theta_s \omega,t) = X(\omega,s+t)$. Consequently, 
$X(\cdot + s)$ has the same law as $X$.  That is, $X$ is stationary.
Once we show that $X$ is  a $(N,d)$
$K$-Brownian motion, this last observation will also apparent from
the observations that, by Theorem 4.9, the random variable $X(\cdot + s)$ is $K$-Gaussian and the support of the law of $X(\cdot + s)$ is the same
as that of $X$.

   Let $\cN$ denote sub-$\sigma$-field of $\cG$ consisting of events with probability either $0$ or $1$.  Set
$$\cF = \sigma\{X_t : t \in D^N\} \vee \cN.$$
Write $\pi_m$ is the quotient map from $D^d$ onto $D^d / \rho^m D^d$
and set
$$\cF^m =  \sigma\{\pi_m \circ X_t : t \in D^N\} \vee \cN.$$
 Similarly, given $C \in \Gamma(N)$, set
$$\cF_C = \sigma\{X_t : t \in C\} \vee \cN$$
and
$$\cF_C^m =  \sigma\{\pi_m \circ X_t : t \in C\} \vee \cN.$$
Also, set
$$\cG_C = \sigma\{Z_B : B \cap C \ne \emptyset\} \vee \cN.$$
Observe that $\cF_C^m \subset \cF_C \subset \cG_C$.

\vskip 24pt
\proclaim Lemma 5.2. Suppose that  $C \in \Gamma(N)$ and $s \in C$.  Then $\cG_C$ and $\sigma\{X_t - X_s : t \notin C\}$ are independent.

\vskip 18pt
\nobreak\noindent
{\it Proof.}  Given $t \notin C$ let $C(m)$ (respectively, $C'(m)$) denote the coset in $\cD(N;m)$ that contains $s$ (respectively, $t$).  Suppose that $|s-t| = q^{-r}$.  Then $C(m)=C'(m)$ for $m \le r$ and $C(m) \cap C'(m) = \emptyset$ for $m>r$.  Moreover, if $m>r$ then, by the isosceles triangle property, $|s-u| = q^{-r}$ for all $u \in C'(m)$ and hence $C'(m) \cap C = \emptyset$.   As $X_t - X_s = \sum_{m=r+1}^\infty [Z_{C'(m)} - Z_{C(m)}]$, it therefore suffices to show that $\cG_C$ and $\{Z_{C'} - Z_{C(m)}: C' \in \cD(N;m), \quad C'\cap C = \emptyset\}$ are independent for all $m$, which will in turn follow if we can show that $Z_{C(m)}$ and   $\{Z_{C'} - Z_{C(m)}: C' \in \cD(N;m), \quad C' \cap C = \emptyset\}$ are independent.  This, however, follows from Corollary 4.12 \hfil\break\rightline{\endproof}

\vskip 24pt
\noindent
{\bf Notation.}  Let $W$ be the $C(D^N, K^d)$-valued random variable  defined by
$$W(\omega)=\cases{\sum_{k=0}^\infty \sum_{C \in \cD(N;k), \, 0 \notin C} Z_C(\omega) \b1_C, &if $\omega\in \Omega^*$, \cr
                     0,                                  &if $\omega\notin \Omega^*$.\cr}$$

\vskip 24pt
\proclaim Lemma 5.3. The $C(D^N, K^d)$-valued random variable
$X - X_0 \b1_{D^N}$  is independent of $X_0$ and has the same law as $W$.

\vskip 18pt
\noindent
{\it Proof.}  Given $C \in \cD(N;n)$ for some $n$, let 
$Z_C' = Z_C - Z_{C(n)}$, where
$C(n) \in \cD(N;n)$ is the ball that contains $0$.  Then
$X - X_0 \b1_{D^N}$ is built from the weights 
$\bigcup_n\{Z_C' : C \in \cD(N;n), \; 0 \notin C\}$
via exactly the same prescription that is used to build $W$
from
$\bigcup_n\{Z_C : C \in \cD(N;n), \; 0 \notin C\}$.
It therefore suffices to show for each $n$ that
$Z_{C(n)}$ is independent of $\{Z_C' :   C \in \cD(N;n), \; 0 \notin C\}$,
and that the latter random variables are i.i.d. with the same common
distribution as the random variables  
$\{Z_C : C \in \cD(N;n), \; 0 \notin C\}$.
This, however, follows from Corollary 4.12.
\hfil\break\rightline{\endproof}

\vskip 24pt
   It is possible to define a rather arbitrary total ordering on $D^N$ as follows.  We begin by inductively defining orders for each of the collections of cosets $\cD(N;n)$.  As $\cD(N;0)$ has only one element, $D^N$, there is certainly no problem in ordering this collection.  Suppose now that $\cD(N;n)$ has been ordered.  For each of the cosets $C$ in $\cD(N;n)$ we assign some arbitrary total ordering to the $q^N$ cosets from $\cD(N;n+1)$ that are contained in $C$.  A total order on $\cD(N;n+1)$ can now be defined such that if $C_1,C_2 \in \cD(N;n+1)$ are contained in distinct cosets in $\cD(N;n)$, say $C_1 \subset C'_1$ and $C_2 \subset C'_2$, then the order relation of $C_1$ and $C_2$ is the same as that of $C'_1$ and $C'_2$.  In
terms of the tree picture, all we are doing is choosing some way of embedding
the tree of balls in the plane and then ordering the balls of a given radius
from left to right.  Using the same embedding, we can also obtain a
total order relation on the boundary of the tree (that is, the points of $D^N$)
by again ordering from left to right.
More formally,
 we declare that the order relation of two distinct points $s$ and $t$ is the same as that for any two disjoint cosets $C_s \ni s$ and $C_t \ni t$ belonging to the same collection $\cD(N;n)$. We will use the symbol $<$ for this ordering on $D^N$.   

\vskip 24pt
\proclaim Lemma 5.4.  Suppose that $t_1,t_2,\ldots,t_n \in D^N$ are such that $t_1 < t_2 < \cdots < t_n$.  Then the random variables $X_{t_1}, X_{t_2} - X_{t_1}, \ldots, X_{t_n} - X_{t_{n-1}}$ are independent.

\vskip 18pt
\nobreak\noindent
{\it Proof.}  It suffices to prove that $X_{t_k} - X_{t_{k-1}}$ is independent of $X_{t_1}, \ldots, X_{t_{k-1}}$ for $2 \le k \le n$.  Suppose that $|t_k - t_{k-1}| = q^{-r}$.  For $1 \le i \le k$ and $m \ge 0$, let $C(i,m)$ denote the coset in $\cD(N;m)$ that contains $t_i$.  Note that $X_{t_k} - X_{t_{k-1}} = \sum_{m=r+1}^\infty [Z_{C(k,m)} - Z_{C(k-1,m)}]$.  From the definition of the ordering, it is clear that $|t_k - t_{k-1}| \le |t_k - t_i|$ for $1 \le i \le k-2$, and hence $C(k,m) \ne C(i,m)$ for $1 \le i \le k-1$ and $m > r$.  The results now follows from Corollary 4.12. \hfil\break\rightline{\endproof}

\vskip 24pt

  We are at last in a position to prove that the law of $X$ has the support $S$,
and so $X$ is a $(N,d)$ $K$-Brownian motion.
Arguing as we did in our informal discussion of the construction of
$X$, we see that $|X_t| \le 1$ for all $t$ and $|X_s - X_t| \le q^{-1}|s-t|$
for all $s,t$.
Thus the support of the distribution of $X$ is contained in $S$.  

  Suppose on the other hand that $f \in S$ and $m \ge 0$.  We will show that
$\bP\{|X_t - f(t)| \le q^{-(m+1)}\} > 0.$
Let $C_1, \ldots, C_{q^{Nm}}$ be a listing of the elements of $\cD(N;m)$ in the order above and choose coset representatives $t_i \in C_i$ for $1 \le i \le q^{Nm}$. As  $|X_s - X_t| \le q^{-1}|s-t|$ and  
$|f(s) - f(t)| \le q^{-1}|s-t|$
for all $s,t$, it follows from
 the ultrametric inequality that $|X_t - f(t)| \le q^{-(m+1)}$ for all $t \in D^N$ if and only if $|X_{t_i} - f(t_i)| \le q^{-(m+1)}$ for $1 \le i \le q^{Nm}$.  By Lemma 3.11, the latter occurs if an only if $|X_{t_1} - f(t_1)| \le q^{-(m+1)}$ and $|X_{t_i} - f(t_i) - (X_{t_{i-1}} - f(t_{i-1}))| \le q^{-(m+1)}$ for $2 \le i \le q^{Nm}$.  Thus, by Lemma 5.4,  
$$\eqalign{&\bP\{ \sup_{t \in D^N} |X_t - f(t)| \le q^{-(m+1)}\}\cr
 &\qquad= \bP\{|X_{t_1} - f(t_1)| \le q^{-(m+1)}\} \prod_{i=2}^{q^{Nm}} \bP\{|X_{t_i} - f(t_i) - (X_{t_{i-1}} - f(t_{i-1}))| \le q^{-(m+1)}\} \cr
   &\qquad= q^{-d(m+1)q^{Nm}} \prod_{i=2}^{q^{Nm}} (q^{-1} |t_i - t_{i-1}|)^{-d} > 0. \cr}$$

\vfill\eject

\centerline{\bf 6. Polar sets, random measures, and $K$-Brownian motion}

\vskip 24pt

   The motivation for much of the material in this section is the question,
``Which Borel sets $B \subset K^d$ are such that 
$\bP\{\exists t \in D^N : X_t \in B \} > 0$.''  The tools that
are used for answering the corresponding question in the 
classical theory of
Brownian motion are basically first hitting times and the strong Markov
propery (see, for example, [Port and Stone, 1978]).  We could define first hitting times in our setting using
one of the orders described prior to Lemma 5.4.  
As Lemma 5.4 shows, $X$ has a nice independent increments
 dependence structure when viewed in
such an order.  However, because our order does not mesh particularly
well with the algebraic structure of $D^N$, this dependence structure 
does not exhibit
any of the time homogeneity that is found in the Brownian case.
Thus, such first hitting times do not appear to be a particularly
suitable tool
for our problem.

   We therefore need some other way of, loosely speaking, picking out
``times'' $t$ such that $X_t \in B$. A similar problem occurs in the
theory of Markov processes with several real parameters
(see, for example, [Dynkin, 1981] and [Fitzsimmons and Salisbury, 1989].)
There the solution is to work with random measures that are supported
on the set of times for which the process is in the set of interest.
Moreover, these random measures are
defined to have suitable ``homogeneity'' and
``adaptedness'' properties.
Random measures and the related additive functionals are, of course, also an important part of the theory of one parameter Markov processes.  Recent works with extensive bibliographies are [Fitzsimmons, 1987] and [Sharpe, 1988].    

   We begin with some definitions.  Let $\cB(D^N)$ denote the Borel $\sigma$-field of $D^N$.  A {\it random measure} is a map $\kappa : \Omega \times \cB(D^N) \rightarrow [0, \infty[$ such that, for each $\omega \in \Omega$, $\kappa(\omega,\cdot)$ is a finite measure on $\cB(D^N)$ and, for each $B \in \cB(D^N)$, $\kappa(\cdot,B)$ is $\cF$-measurable.  As usual, we will often write $\kappa(B)$ for $\kappa(\cdot,B)$.  

   A random measure $\kappa$ is {\it integrable} if
$$\bP[\kappa(D^N)] < \infty$$ 
and {\it square-integrable} if 
$$\bP[\kappa^2(D^N)] < \infty.$$
Given an integrable random measure $\kappa$, there is a finite Borel measure $\mu$ on $K^d$ defined by $\int f(x) \,\mu(dx) = \bP [\int f(X_t)\, \kappa(dt)]$.   Very loosely speaking, if we pick a point $t \in D^N$ ``at random'' according to $\kappa$ (of course, $\kappa$ is not necessarily
a probability measure), then the ``distribution'' of $X_t$ is $\mu$
(again, $\mu$ is not necessarily a probability measure). We say that $\mu$ is the {\it characteristic measure} of $\kappa$.  We will see that $\mu$ plays
a role similar to  that of the Revuz measure does in the 
Markov process theory.  

   A random measure $\kappa$ is {\it adapted} if $\kappa(C)$ is $\cF_C$-measurable for all $C \in \Gamma(N)$. (Recall the $\Gamma(N)$
is the collection of sub-balls of $D^N$.) 

   Given a random measure $\kappa$ and $t \in D^N$ we can define another random measure $\Theta_t \kappa$ by setting $\Theta_t \kappa(\omega,B) = \kappa(\theta_t \omega, t+B)$.  A random measure $\kappa$ is {\it homogeneous} if $\Theta_t \kappa = \kappa$ for all $t \in D^N$.

   The most obvious examples of an integrable, adapted, homogeneous random
measure are the measures of the form $\kappa(B) = \int_B h(X_t), dt$ for 
some bounded, nonnegative, Borel function $h$.  As each random variable
$X_t$ is $K$-Gaussian with support $D^d$, it follows
by Fubini's theorem that $\kappa$ has characteristic measure
$\mu(dx) = h(x) \b1_{D^d} (x) dx$.  We will show in Theorem 6.1 that
all integrable, adapted, homogeneous random measures arise as limits
of random measures of this type.

   When we say that a particular random measure $\kappa$ is the unique random measure with  certain properties, we mean that if $\kappa'$ is another random measure with the same properties then the two measure $\kappa(\omega,\cdot)$ and $\kappa'(\omega,\cdot)$ coincide for almost all $\omega \in \Omega$.

\vskip 24pt

\noindent
{\bf Notation.} Define an approximate identity 
$\{\phi_n : K^d \rightarrow \bR\}_{n=0}^\infty$ by setting
$$\phi_n(x) = \cases{q^{dn}, & if $|x| \le q^{-n}$, \cr
                     0,      & otherwise. \cr}$$
Thus $\phi_n(x) dx$ is the $K$-Gaussian measure with support $\rho^n D^d$.

\vskip 24pt
\proclaim Theorem 6.1.  Suppose that $\kappa$ is an integrable, adapted, homogeneous random measure with characteristic measure $\mu$.  Then, for all $B$ in the algebra consisting of unions of sub-balls of $D^N$
(that is, the algebra generated by $\Gamma(N)$), the sequence of random variables $\{\int_B \phi_n \ast \mu(X_t)\, dt\}_{n=0}^\infty$ is a martingale 
for the filtration $\{\cF^n\}_{n=0}^\infty$.  Each such martingale
converges almost surely and in $\cL^1$ to
$\kappa(B)$ as $n\rightarrow\infty$.  In particular, $\kappa$ is the unique integrable, adapted, homogeneous random measure with characteristic measure $\mu$.
  
\vskip 18pt
\nobreak\noindent
{\it Proof.}  Fix an integer $n \ge 0$.
Recall from \S 5 that $\cD(N;n)$ is the collection of cosets of
$\rho^n D^N$ in $D^N$.  
Observe from the construction of $X$ that if $C$ is a coset belonging to $\cD(N;n)$ then the random function $\pi_n \circ X$ 
(where we remind the reader that $\pi_n$ is the quotient map
from $D^d$ onto $D^d/\rho^n D^d$) is constant on $C$.  Denote the common value of $\pi_n \circ X_t$, $t \in C$, by $Y_C$.  As the law of $X_t$
is normalised Haar measure on $D^d$, the law of $Y_C$ is just normalised Haar measure on $D^d / \rho^n D^d$.  As $\kappa$ is homogeneous, there is some function $g_n : D^d / \rho^n D^d \rightarrow [0,\infty[$ such that $q^{Nn}\bP[\kappa(C)|\cF_C^n] = g_n(Y_C)$ for all such $C$.  Given  any $f : D^d / \rho^n D^d \rightarrow \bR$  we have
$$\eqalign{\int f \circ \pi_n(x) \,\mu(dx) &= \bP [\int f \circ \pi_n(X_t)\, \kappa(dt)] \cr
                        &= \sum_C \bP[f(Y_C) g_n(Y_C)] \cr
                        &= q^{-dn} \sum_{y \in D^d / \rho^n D^d} f(y) g_n(y). \cr}$$
Thus $g_n(y) = q^{dn} \mu(\pi_n^{-1}y)$.  Equivalently, $g_n(\pi_n x) = \phi_n \ast \mu (x)$ for all $x \in D^d$.

   From the homogeneity of $\kappa$ it follows that if $C' \in \cD(N;m)$ for some $m \ge n$ and $C' \subset C \in \cD(N;n)$ then 
$$\bP[\kappa(C') | \cF_C^n] = q^{-N(m-n)} \bP[\kappa(C) | \cF_C^n] = q^{-Nm}g_n(Y_C).$$  
Finally, note that for such a pair of cosets $C$ and $C'$ we have  $\cF^n = \cF_C^n \vee \sigma\{\pi_n \circ X_t - \pi_n \circ X_s : t \notin C\}$, for any fixed $s \in C$; and  hence, from Lemma 5.2 and the adaptedness of $\kappa$, we have
$$\bP[\kappa(C') | \cF_C^n] = \bP[\kappa(C') | \cF^n].$$

   Combining all of the above observations, we see that for any $B$ in the algebra generated by $\Gamma(N)$ we have
$$\bP[\kappa(B) | \cF^n] = \int_B \phi_n \ast \mu(X_t)\, dt.$$
As $\bigvee_{n=0}^\infty \cF^n = \cF$, the result now follows from  the martingale  convergence theorem.  \hfil\break\rightline{\endproof}

\vskip 24pt

   The following {\it Palm  formula} is our counterpart of the strong
Markov property. Recall from Lemma 5.3 that $X - X_0 \b1_{D^N}$ is
independent of $X_0$ and has the same law as $W$.  Consequently,
for a fixed $t \in D^N$,
$X - X_t \b1_{D^N}$ is independent of $X_t$ 
and has the same law as $W(\cdot - t)$. In the same very loose manner
that we interpreted the meaning of the characteristic measure of
an integrable random measure, the Palm formula can be thought of
as saying that if we pick a point
$t \in D^N$
``at random'' according to an integrable, adapted, homogeneous random
measure $\kappa$ with characteristic measure $\mu$, then $X_t$
has ``distribution'' $\mu$ and the ``conditional law'' of 
$X - X_t \b1_{D^N}$ given $X_t$ is that of $W(\cdot -t)$.

\vskip 24pt
\proclaim Corollary 6.2.  Suppose that $\kappa$ is an integrable, adapted, homogeneous random measure with characteristic measure $\mu$, and $F$ is a non-negative, Borel function on the space of continuous functions $C(D^N,K^d)$.  Then
$$\bP [\int F(X \circ \theta_t) \,\kappa(dt)] = \int \bP[F(x+W)] \,\mu(dx).$$

\vskip 18pt
\nobreak\noindent
{\it Proof.}  By a monotone class argument, it suffices to prove the result when $F$ is of the form $F(f) = G(f(t_1),\ldots,f(t_k))$ for some finite set of points $t_1, \ldots, t_k \in D^N$ and some bounded, continuous function $G:(K^d)^k \rightarrow \bR$.

   From Theorem 6.1, we know that the sequence of random measure $\{\kappa_n\}_{n=0}^\infty$ defined by $\kappa_n(B) = \int_B \phi_n \ast \mu(X_t)\, dt$ converges weakly to $\kappa$ almost surely.  Moreover, $\kappa_n(D^N)$ converges to $\kappa(D^N)$ in $\cL^1$ as $n \rightarrow \infty$.  From Lemma 5.3 we see that
$$\eqalign{\bP[\int F(X \circ \theta_t) \,\kappa_n(dt)]
            &= \bP[\int G(X_{t+t_1},\ldots,X_{t+t_k})\, \phi_n \ast \mu(X_t)\, dt] \cr
            &= \bP[G(X_{t_1},\ldots,X_{t_k}) \,\phi_n \ast \mu(X_0)] \cr
            &= \int \bP[G(x+W_{t_1},\ldots,x+W_{t_k})]\, \phi_n \ast \mu(x)\, dx \cr
            &= \int \bP[F(x+W)] \,\phi_n \ast \mu(x) \,dx. \cr}$$
The result follows when we let $n \rightarrow \infty$ and appeal to dominated convergence.  \hfil\break\rightline{\endproof}

\vskip 24pt

\vskip 24pt
\noindent
{\bf Notation.}  Define a function $u:\{q^{-k}\}_{k=0}^\infty \cup \{0\} \rightarrow [0,\infty]$ as follows.  For $r \ne 0$ set
$$u(r) = \cases{{{1-q^{-N}} \over {q^d}} \log_q {1 \over r}, & if $N=d$, \cr
                 {{1-q^{-N}} \over {q^d}} {1 \over {q^{(d-N)}-1}} [{1 \over {r^{(d-N)}}}-1], & if $N \ne d$. \cr}$$
Set
$$u(0) = \cases{{{1-q^{-N}} \over {q^d}} {1 \over {1-q^{(d-N)}}}, & if $N>d$, \cr
                \infty,                                         & if $N \le d$.
 \cr}$$

\vskip 24pt

The functions $(x,y) \mapsto u(|x-y|)$ appear in the local field
Riesz potential theory of [Taibleson, 1975] (see, also, [Evans, 1988b, 1992]).
They will turn out to play the role
of the object variously referred to in  Markov process theory as the
$0$-potential kernel density, the $0$-resolvent density,  or the Green's function.

The second part of the following  theorem is the fundamental existence result for square-integrable, adapted, homogeneous random measures.  The corresponding result in the integrable case is given in Theorem 6.10.

\vskip 24pt
\proclaim Theorem 6.3.  Suppose that $\kappa$ is a square-integrable, adapted, homogeneous random measure with characteristic measure $\mu$.  Then
$$\bP[\kappa^2(D^N)] = \int \int u(|x-y|) \,\mu(dx) \mu(dy).$$
Conversely, given any finite measure $\mu$ on $D^d$ such that $\int \int u(|x-y|) \,\mu(dx) \mu(dy) < \infty$, there exists a unique,  square-integrable, adapted, homogeneous random measure $\kappa$ with characteristic measure $\mu$.

\vskip 18pt
\nobreak\noindent
{\it Proof.} Suppose first of all that $\kappa$ is a square-integrable, adapted, homogeneous random measure with characteristic measure $\mu$.  It follows from Theorem 6.1 and the $\cL^2$ martingale convergence theorem that
$$\eqalign{\bP[\kappa^2(D^N)] &= \lim_{n \rightarrow \infty} \int \int \bP[\phi_n \ast \mu(X_s) \, \phi_n \ast \mu(X_t)]\, ds dt \cr
                              &= \lim_{n \rightarrow \infty} \int \bP[\phi_n \ast \mu(X_0) \, \phi_n \ast \mu(X_t)]\, dt. \cr}$$
We know from Lemma 5.3 that $X_0$ and $X_t - X_0$ are independent.
The support of the law of $X_0$ (resp. $X_t - X_0$) is $D^d$ (resp.
$\{z \in D^d : |z| \le q^{-1} |t|\}$).  The support of the law of $(X_0, X_t)$ is thus the
$D$-submodule 
$$\{(w,w+z) \in (K^d)^2 : |w| \le 1, \; |z| \le q^{-1} |t|\}
= \{(x,y) \in (K^d)^2 : |x| \le 1, \; |x-y| \le q^{-1} |t|\}.$$
A simple integration shows that Haar measure on $(K^d)^2$ assigns mass
$(q^{-1} |t|)^d$ to this set, and so 
 the law of the $K$-Gaussian pair $(X_0, X_t)$ 
has a density with respect to Haar measure on $(K^d)^2$ given by 
$${{\bP\{X_0 \in dx, X_t \in dy\}} \over {dx \, dy}} = \cases{(q^{-1}|t|)^{-d}, & if $|x| \le 1$, $|y| \le 1$ and $|x-y| \le q^{-1}|t|$, \cr
    0, & otherwise. \cr}$$
Another integration shows that
$$\eqalign{\int \bP[\phi_n \ast \mu(X_0) \, \phi_n \ast \mu(X_t)]\, dt &= \int \int u(|x-y|) \phi_n \ast \mu(x) \, \phi_n \ast \mu(y) \, dx dy \cr
   &= \int \int \phi_n \ast v(x-y) \, \mu(dx) \mu(dy), \cr}$$
where we set $v(z) = u(|z|)$.  Note that if $|z|=q^{-m}$ then $\phi_n \ast v(z) = v(z) = u(|z|)$ for all $n > m$.  Also, $\phi_n \ast v(0) \uparrow v(0) = u(|0|)$ as $n \rightarrow \infty$.  Thus 
$$\lim_{n \rightarrow \infty} \int \int \phi_n \ast v(x-y) \, \mu(dx) \mu(dy) =  \int \int u(|x-y|) \, \mu(dx) \mu(dy),$$
as required.

   Conversely,  suppose that $\mu$ is a finite measure such that $\int \int u(|x-y|)\, \mu(dx) \mu(dy) < \infty$.  Define a sequence $\{\kappa_n\}_{n=0}^\infty$ of square-integrable, adapted, homogeneous random measure by setting
$$\kappa_n(A) = \int_A \phi_n \ast \mu(X_t)\, dt.$$
By an argument similar to that given in the proof of Theorem 6.1, $\{\kappa_n(B)\}_{n=0}^\infty$ is a martingale for each $B$ in the algebra generated by $\Gamma(N)$.  Moreover, from the calculations of the previous paragraph it is clear that each of these martingales is bounded in $\cL^2$ and hence convergent almost surely and in $\cL^2$.  Let $\Omega'$ denote the subset of $\Omega$ of probability one on which all of these martingales converge.  Define a set function $\kappa$ on  the algebra generated by $\Gamma(N)$ by setting
$$\kappa(\omega,B) = \cases{\lim_{n \rightarrow \infty} \kappa_n(\omega,B), & if $\omega \in \Omega'$, \cr
                            0, & otherwise. \cr}$$

   Note that $\kappa(\omega,\cdot)$ is finitely additive for all $\omega \in \Omega$.  Note also that if $B_0 \supset B_1 \supset \ldots$ is a decreasing sequence of sets in this algebra such that $\cap_i B_i =\emptyset$ then we certainly have  $\lim_{i \rightarrow \infty} \kappa(\omega, B_i) = 0$ for all $\omega \in \Omega$, as each set $B_i$ is compact and hence $B_i = \emptyset$ for all sufficiently large $i$.  Applying a standard extension theorem (see, for example, Theorems 3.1.1 and 3.1.4 of [Dudley, 1989]), we see that $\kappa(\omega, \cdot)$ extends uniquely to a finite measure on $\cB(D^N)$ for all $\omega \in \Omega$.  A monotone class argument shows that $\kappa$ is an adapted random measure.  As each of the random measures $\kappa_n$  is homogeneous and $\theta_t \Omega' = \Omega'$ for all $t$, it is clear that $\kappa$ is homogeneous.  The algebra generated by $\Gamma(N)$ is a weak convergence determining class
(that is, if $\{\nu_n\}_{n=0}^\infty$ is a sequence of finite measures such
that $\nu_n(B)$ converges as $n \rightarrow \infty$ to $\nu(B)$ for some
finite measure $\nu$, then $\nu_n$ converges weakly to $\nu$ as $n \rightarrow
\infty$, cf. \S 3.4 of [Ethier and Kurtz, 1986]).
Hence $\kappa_n(\omega,\cdot)$ converges weakly to $\kappa(\omega,\cdot)$ for almost all $\omega \in \Omega$.  Thus, for any bounded, continuous function $f:K^d \rightarrow \bR$ we have that $\int f(X_t) \,\kappa(dt) = \lim_{n \rightarrow \infty} \int f(X_t) \,\kappa_n(dt)$ almost surely and hence, by dominated convergence,
$$ \eqalign{\bP[\int f(X_t) \,\kappa(dt)] &= \lim_{n \rightarrow \infty} \bP[\int f(X_t) \,\kappa_n(dt)] \cr
                                        &= \lim_{n \rightarrow \infty} \int f(x) \, \phi_n \ast \mu(x) \, dx \cr
                                        &= \int f(x)\, \mu(dx). \cr}$$
Another monotone class argument then shows that $\kappa$ has characteristic measure $\mu$.  Applying Theorem 6.1, we see that $\kappa$ is the unique such adapted, homogeneous random measure.  \hfil\break\rightline{\endproof}

\vskip 24pt

   Recall that we want to answer the question posed at the beginning of the
section, and we have introduced random measures as a technique for
somehow picking out ``times'' $t \in D^N$ such that $X_t \in B$.
Of course, it is not hard to build such a random measure - one can
just take the unit point mass at the ``first hitting time'' of $B$ in one of
the orders introduced in Section \S 5. (This will work if $B$ is
closed).  However, as we have observed, such a random measure is not
particularly consonant with the structure of $X$.  In order to
arrive at an answer to our question, we will need to transform such a
first hitting time into a square-integrable, adapted, homogeneous,
random measure - the sort of measure for which we have a good
analogue of the strong Markov property (recall Corollary 6.2).
We will carry out such a transformation in the proof of Theorem 6.5,
where we first turn the unit point mass at
our hitting time into a homogeneous random measure by a simple
averaging procedure and then 
use the following  analogue of the `central projection' operation of [Dynkin, 1981] that converts non-adapted random measures into adapted ones.

\proclaim Theorem 6.4.  Suppose that $\kappa$ is a square-integrable, homogeneous random measure with characteristic measure $\mu$.  Then there exists a unique square-integrable, adapted, homogeneous random measure $\tilde \kappa$ with characteristic measure $\mu$.

\vskip 18pt
\nobreak\noindent
{\it Proof.}  By Theorem 6.3, it will suffice to show that $\int \int u(|x-y|) \, \mu(dx) \mu(dy) < \infty$.  Given an integer $n \ge 0$,  define an integrable, homogeneous random measure $\kappa_n$  by setting 
$$\kappa_n(dt) = \bP[\kappa(C)|\cF_C^n] q^{Nn} \lambda^N(dt),$$ 
where $C \in \cD(N;n)$ is the coset that contains $t$.  By the homogeneity of $\kappa$ we see that there is a function $h_n : D^d / \rho^n D^d \rightarrow \infty$ such that for all $C \in \cD(N;n)$ and any $t \in C$ we have
$$\bP[\kappa(C)|\cF_C^n] = h_n(\pi_n \circ X_t),$$
and hence $\kappa_n$ is adapted.  An argument similar to that in the proof of Theorem 6.1 gives that
$$h_n(y) = q^{(d-N)n} \mu(\pi_n^{-1} y)$$
and hence $\kappa_n$ has characteristic measure $\phi_n \ast \mu$.  From Theorem 6.3 we have
$$\bP[\kappa_n^2(D^N)] = \int \int \phi_n \ast v(x-y) \, \mu(dx) \mu(dy),$$
where we set $v(z) = u(|z|)$.  The required result will follow from Fatou's lemma if we can show that $\bP[\kappa_n^2(D^N)]$ is uniformly bounded.

   Let $C_1,\ldots,C_{q^{Nn}}$ be a listing of the cosets in $\cD(N;n)$ in 
one of the orders described in \S 5.  Set
$$\cI_i = \bigvee_{j=1}^i \cG_{C_j}$$
and    
$$\cJ_i = \bigvee_{j=i}^{q^{Nn}} \cG_{C_j}.$$
Observe that
$$\bP[\bP[\kappa(C_i)|\cI_i] | \cJ_i] = \bP[\kappa(C_i) | \cG_{C_i}],$$
(cf. Lemma A.4 of the Appendix.)  Thus, by two applications of a discrete time version of Meyer's energy inequality (cf. Lemma A.5 of the Appendix), we have 
$$\bP[\{\sum_i \bP[\kappa(C_i) | \cG_{C_i}]\}^2] \le 16 \bP[\kappa^2(D^N)].$$
Furthermore, if $s \in C_i$ then $\cF = \cF_{C_i} \vee \sigma\{X_t-X_s : t \notin C_i\}$, and it follows from Lemma 5.2 that $\kappa_n(C_i) = \bP[\bP[\kappa(C_i) | \cG_{C_i}] | \cF]$, and hence, by Jensen's inequality for conditional expectations, $\bP[\kappa_n^2(D^N)] \le 16\bP[\kappa^2(D^N)]$.  \hfil\break\rightline{\endproof}

\vskip 24pt
   We have now marshalled enough tools to answer the question posed at the beginning of the section.  For convenience we will say, as usual, that a Borel subset $B \subset D^d$ is {\it polar} if $\bP\{\exists t \in D^N : X_t \in B\} = 0$.  Otherwise, $B$ is {\it non-polar}.  The following necessary
and sufficient condition for polarity is analogous to the well-known results of [Kakutani, 1944a, 1944b] in the Brownian case.

\vskip 24pt
\proclaim Theorem 6.5.  A Borel set $B \subset D^d$ is non-polar if and only if there is a non-trivial finite measure $\mu$ such that $\mu(D^d \backslash B) = 0$ and $\int \int u(|x-y|) \,\mu(dx) \mu(dy) < \infty$.

\vskip 18pt
\nobreak\noindent
{\it Proof.}  Suppose that $B \subset D^d$ is non-polar.  The map
$A \mapsto \bP\{\exists t \in D^N : X_t \in A\}$ defined
on the Borel sets of $D^d$ is a Choquet capacity
relative to the compact paving, and Choquet's capacitability theorem 
gives that
$$\bP\{\exists t \in D^N : X_t \in A\}
 = \sup_K \bP\{\exists t \in D^N : X_t \in K\}, $$
where the supremum is over compact subsets of $A$
(cf. Appendix II.3 of [Doob, 1984]).  
Thus there is a compact subset of $B$ which is also non-polar, so we may suppose without loss of generality that $B$ is compact.

   Let $(D^N)^\Delta=D^N\cup\{\Delta\}$, where $\Delta$ is an abstract isolated point.  Using one of the orders on $D^N$ introduced in \S 5, define a $\cF$-measurable map $S:\Omega\rightarrow (D^N)^\Delta$ 
by setting
$$S(\omega) = \cases{ \min\{t \in D^N: X_t \in B\}, &if $X(\omega,D^N) \cap B \ne \emptyset$, \cr
\Delta, &otherwise. \cr}$$ 
Thus
$$(\omega,S(\omega))\in\{(\omega',t):X(\omega',t)\in B\} \Leftrightarrow X(\omega,D^N)\cap B\ne\emptyset $$
and
$$S(\omega)=\Delta \Leftrightarrow X(\omega,D^N)\cap B=\emptyset.$$
  Given $\omega\in\Omega$ and a Borel set $A \subset D^N$, define $\kappa(\omega,A)=\int {\bf 1}_A((S\circ\theta_t)(\omega)+t) \,dt$ for $\omega$ such that $X(\omega,D^N)\cap B\ne\emptyset$ and put $\kappa(\omega,A)=0$ otherwise.  It is straightforward to check that $\kappa$ is a homogeneous random measure with the following properties:
$$\bP[\kappa(D^N)]=\bP\{\exists t\in D^N:X_t\in B\},$$
$$\kappa(\{t\in D^N : X_t\notin B\})=0,$$
$$\kappa(D^N) \le 1.$$
In particular, $\kappa$ is square-integrable and has a non-trivial characteristic measure $\mu$ that is supported on $B$.  Applying Theorems 6.4 and 6.3, we see that   $\int \int u(|x-y|) \,\mu(dx) \mu(dy) < \infty$.

   Conversely, suppose that a Borel set $B \subset D^d$ is such that there is a non-trivial finite measure $\mu$ satisfying $\mu(D^d \backslash B) = 0$ and $\int \int u(|x-y|) \,\mu(dx) \mu(dy) < \infty$.   From Theorem 6.3 we know that there is a square-integrable, adapted, homogeneous random measure with characteristic measure $\mu$, and this certainly implies that $B$ is non-polar. \hfil\break\rightline{\endproof}

\vskip 24pt

   The only measures supported by a singleton are, of course, the
corresponding
point masses, and so Theorem 6.5 gives that singletons will be non-polar
if and only if $u(0) < \infty$ -  that is, if and only if $N > d$.   
When $N \le d$ we have the following criteria for polarity and non-polarity.

\vskip 24pt
\noindent
{\bf Notation.}  Given a non-decreasing function $h:\{q^{-n}\}_{n=0}^\infty \rightarrow [0,\infty[$ with $\lim_{n \rightarrow \infty} h(q^{-n}) = 0$,
let $h-m(\cdot)$ denote the {\it Hausdorff measure} on $K^d$ constructed using the measure function $h$.  That is, for a Borel set $A \subset K^d$,
$$h-m(A) = \lim_{n \rightarrow \infty} \inf\{\sum_i h({\rm diam} B_i)\},$$
where the infimum is taken over all countable collections of balls $\{B_i\}$
such that $A \subset \bigcup_i B_i$ and $\sup_i {\rm diam} B_i \le q^{-n}$.
When $h(r) = r^\alpha$ for $\alpha > 0$, we will write $r^\alpha-m$ for $h-m$.
The {\it Hausdorff dimension} of a set A is given by
$$\inf\{\alpha : r^\alpha-m(A) = 0\}=\sup\{\alpha : r^\alpha-m(A) = \infty\}.$$
(It is not hard to show that the measure $r^d-m$ coincides with the trace of $\lambda$ on $D^d$, and
$D^d$ has Hausdorff dimension $d$.) 

\vskip 24pt
\proclaim Corollary 6.6.  Suppose that $N \le d$.  If $B$ is a Borel subset of $D^d$ with Hausdorff dimension less than $d-N$ then $B$ is polar.  On the other hand, if $B$ has Hausdorff dimension greater than $d-N$ then $B$ is non-polar.

\vskip 18pt
\nobreak\noindent
{\it Proof.} Given Theorem 6.5, the result follows immediately from a local
field version of Frostman's theorem connecting Riesz-type capacities
and Hausdorff dimension (see Theorem 2.3 of [Evans, 1988b] for more details.)  \hfil\break\rightline{\endproof}

\vskip 24pt

   It is clear that the characteristic measure of an integrable, adapted, homogeneous random measure is finite and does not charge any polar sets.  We 
will finish this section by showing in Theorem 6.10 below that the converse holds.
We first need the following three analytic results that are similar ones
found in classical theory of Riesz potentials (cf. \S\S I.3 and I.4 of [Landkof, 1972]), although our proofs are
somewhat different 

\vskip 24pt
\proclaim Lemma 6.7.  Suppose that $\mu$ is a finite measure on $D^N$ with support $A$.  Then
$$\sup\{\int u(|x-y|) \mu(dy) : x \in D^N\} = \sup\{\int u(|x-y|) \mu(dy) : x \in A\}.$$

\vskip 18pt
\nobreak\noindent
{\it Proof.} Given $x \notin A$, we may choose $z \in A$ such that $|x-z| = \inf\{|x-y| : y \in A\}$.  For any $y \in A$ we have $|x-z| \le |x-y|$ and hence, by the strong triangle inequality, $|z-y| \le \max\{|x-z|, |x-y|\} \le |x-y|$.  Thus $u(|z-y|) \ge u(|x-y|)$ for all $y \in A$, and $\int u(|z-y|)\, \mu(dy) \ge \int u(|x-y|)\, \mu(dy)$, as required.  \hfil\break\rightline{\endproof}

\vskip 24pt
\proclaim Lemma 6.8.  Suppose that $\mu$ is a non-trivial, finite measure on $D^N$ such that $\mu(\{x : \int u(|x-y|)\, \mu(dy) = \infty\}) = 0$.  Then there exists a non-trivial measure $\nu \le \mu$ such that the function $x \mapsto \int u(|x-y|) \,\nu(dy)$ is uniformly bounded.

\vskip 18pt
\nobreak\noindent
{\it Proof.} By assumption, we may choose a compact set $B$ such that $\mu(B) > 0$ and 
$$\sup\{\int u(|x-y|)\, \mu(dy) : x \in B\} < \infty.$$
Now let $\nu$ be the trace of $\mu$ on $B$ and use Lemma 6.7.  \hfil\break\rightline{\endproof}

\vskip 24pt
\proclaim Lemma 6.9.  Suppose that $\mu$ is a finite measure on $D^N$ that does not charge any polar set.  Then
$$\mu = \sup\{\nu \le \mu : \int \int u(|x-y|) \,\nu(dx) \nu(dy) < \infty\}.$$

\vskip 18pt
\nobreak\noindent
{\it Proof.} Let $\sigma = \sup\{\nu \le \mu : \int \int u(|x-y|)\, \nu(dx) \nu(dy) < \infty\}$ and suppose that $\sigma \ne \mu$.  Let $\eta = \mu - \sigma$.

   Note first of all that $B = \{x : \int u(|x-y|) \,\eta(dx) = \infty\}$ is polar.  For if this was not the case then, by Theorem 6.5, there would exist a non-trivial, finite measure $\xi$ concentrated on $B$ such that $\int \int u(|x-y|) \,\xi(dx) \xi(dy) < \infty$.  Applying Lemma 6.8, we may even suppose that $x \mapsto \int u(|x-y|)\, \xi(dy)$ is uniformly bounded, say by $c$. Then, however, we would have the contradiction
$$\eqalign{ c \eta(D^N) &\ge \int \int u(|x-y|) \,\xi(dy) \eta(dx)\cr
                        &= \int \int u(|x-y|)\, \eta(dx) \xi(dy) = \infty.\cr}$$

   As $B$ is polar and $\eta \le \mu$, we have $\eta(B) = 0$.  Another application of Lemma 6.8 shows that there exists a non-trivial, finite measure $\nu \le \eta$ such that $ x \mapsto \int u(|x-y|) \,\nu(dy)$ is uniformly bounded, and hence $\int \int u(|x-y|)\, \nu(dx) \nu(dy) < \infty$.

   Finally, suppose that $\nu_1 \le \nu_2 \le \cdots \le \mu$ is such that $\int \int u(|x-y|) \,\nu_n(dx) \nu_n(dy) < \infty$ for all $n$ and $\nu_n \uparrow \sigma$ as $n \rightarrow \infty$.  Then $\nu_n + \nu \le \mu$ for all $n$ and $\int \int u(|x-y|)\, (\nu_n + \nu)(dx) (\nu_n + \nu)(dy) < \infty$ (cf. the Cauchy-Schwarz inequality in \S 3 of [Evans, 1988b].)  However, $\nu_n + \nu \uparrow \sigma + \nu > \sigma$, which contradicts the definition of $\sigma$.  \hfil\break\rightline{\endproof}

\vskip 24pt
\proclaim Theorem 6.10.  Suppose that $\mu$ is a finite measure on $D^N$ that does not charge any polar set.  Then $\mu$ is the characteristic measure of a unique integrable, adapted, homogeneous random measure.

\vskip 18pt
\nobreak\noindent
{\it Proof.}  From Lemma 6.9 we know that there is a sequence of finite measures $\mu_1 \le \mu_2 \le \cdots \le \mu$ such that $ \int \int u(|x-y|) \,\mu_n(dx) \mu_n(dy) < \infty $ for all $n$ and $\mu_n \uparrow \mu$ as $n \rightarrow \infty$.  From Theorem 6.4, we have that each measure $\mu_n$ is the characteristic measure of a square-integrable, adapted, homogeneous random measure $\kappa_n$.  Furthermore, each  measure $\mu_{n+1} - \mu_n$ is also  the characteristic measure of a square-integrable, adapted, homogeneous random measure, and so we see from Theorem 6.1 that we may suppose that $\kappa_1(\omega,\cdot) \le \kappa_2(\omega,\cdot) \le \cdots$ for all $\omega \in \Omega$.  Set $\kappa = \sup_n \kappa_n$.  By the monotone convergence theorem, $\kappa$ is an integrable random measure with characteristic measure $\mu$.  Moreover, it is clear that $\kappa$ is adapted and homogeneous.  The uniqueness of $\kappa$ follows from Theorem 6.1.  \hfil\break\rightline{\endproof}

\vfill\eject
\centerline{\bf 7. Constructions and continuity of $K$-Brownian local time}

\vskip 24pt

   For the whole of this section we will continue to use the notation $X$
for the  particular $(N,d)$ $K$-Brownian motion constructed in \S 5.
Assume that $N > d$, and hence the singleton subsets of $D^d$ are all non-polar.  Given $x \in D^d$,  write $\delta_x$ for the unit point mass at $x$.  From Theorem 6.3 we see that there exists a unique square-integrable, adapted, homogeneous random measure $L_x$ with characteristic measure $\delta_x$.  In keeping with the Markov process nomenclature, we say that $L_x$ is the {\it local time at the level} $x$.  Recall from the general construction
in the proof of Theorem 6.1 that the random 
measure $L_x$ is the almost sure weak limit of the sequence of random measures $\{L_x^{(n)}\}_{n=0}^\infty$ defined by $L_x^{(n)}(B) = \int_B \phi_n(X_t - x) \, dt$.  As the function $x \b1_{D^N}$ belongs to the support
of $X$, it follows that $X-x \b1_{D^N}$ has the same law as $X$, and hence $L_x$ has the same distribution as $L_0$.  In discussing the properties of local times at a fixed level $x$ it therefore suffices to consider the case $x=0$.  

   The local time $L_0$ is almost surely supported on the zero set $\cZ=\{t \in D^N : X_t = 0\}$ (in fact, by Corollary 7.3 below the closed support of $L_0$ is $\cZ$), and it is of interest to investigate whether there are alternative constructions of $L_0$ that obtain $L_0$ by applying some deterministic operation to $\cZ$.  In the Markov process literature such constructions are described as being {\it intrinsic}.  Our intuition is that such a construction
should proceed by spreading mass as smoothly and as evenly on $\cZ$ as possible.

   If $\cZ$ was a set of positive Haar measure, then the obvious candidate
for such a procedure would be to simply use the trace of Haar measure on
$\cZ$.  However, since $\bP\{X_t = 0\} = 0$ for all $t \in D^N$, it follows
from Fubini's theorem that $\cZ$ has Haar measure $0$ almost surely, 
and this idea
fails.  We can salvage some of the features of this procedure in two different
ways.  The first approach is to  build the $\epsilon$-neighbourhood 
of $\cZ$, take the trace of Haar measure on this latter
set, let
$\epsilon \downarrow 0$, and hope to be able to renormalise the resulting
sequence of measures so that we get $L_0$ in the limit.  This is the analogue
of the dilation construction of Brownian local time in [Kingman, 1973],
and we carry it out in Theorem 7.2.  The second approach is to use
the trace on $\cZ$ of a measure that is translation invariant like
Haar measure, but assigns finite, positive mass to $\cZ$.
This is the analogue of the Hausdorff measure construction of
[Taylor and Wendel, 1966] for Brownian local time.
We describe this result without proof in Theorem 7.4. A proof
may be found in [Evans, 1993].

  We do not consider the more delicate question of whether these constructions produce the local time for all levels simultaneously  (see, for example, [Barlow et al., 1986] for a discussion of this question in a Euclidean setting.)

  In Theorem 7.5 we turn to the study of $L$ as a measure-valued process
indexed by the level at which the local time is evaluated, and 
obtain an analogue of the theorem of [Trotter, 1958] that the local
time of linear Brownian motion is continuous in space and time.  

The following observation is fundamental.

\vskip 24pt
\proclaim Lemma 7.1.  Suppose that $\kappa$ is an integrable, adapted, homogeneous random measure that is almost surely supported on $\cZ$.  Then $\kappa = cL_0$ for some constant $0 \le c < \infty$.

\vskip 18pt
\nobreak\noindent
{\it Proof.}  By assumption, the characteristic measure of $\kappa$ is of the form $c \delta_0$ for some constant  $0 \le c < \infty$, and the result follows from Theorem 6.1. \hfil\break\rightline{\endproof}

\vskip 24pt

\noindent
{\bf Notation.} For an integer $n \ge 0$ let $\cZ^{(n)}$ be the 
$q^{-n}$-neighbourhood of $\cZ$.  That is,
$$\cZ^{(n)} = \rho^n D^N + \cZ = \{t \in D^N : |t-\cZ| \le q^{-n}\}.$$
Note that $\cZ^{(n)}$ is the union of all the balls in $\cD(N;n)$ that intersect
$\cZ$. 

\vskip 24pt

\proclaim Theorem 7.2.  The sequence of random measures $\{\kappa_n\}_{n=0}^\infty$ defined by 
$$\kappa_n(B) = q^{dn} \lambda(\cZ^{(n)} \cap B)$$
converges weakly  to $hL_0$ almost surely, where $h = \bP[0 \in X(D^N)] = \bP[\cZ \ne \emptyset]$.

\vskip 18pt
\nobreak\noindent
{\it Proof.}   For an integer $n \ge 0$ let $V_n$ denote the number of cosets in $D(N;n)$ that intersect $\cZ$.  Note that $\kappa_n(D^N) = q^{(d-N)n} V_n$. 

   Given a coset $C \in \cD(N;1)$ define a $C(C, K^d)$-valued
random variable $X^C$ by setting $X_t^C = X_t - Z_{D^N} \b1_C$ for $t \in C$.
In the notation of \S 5,
$$X^C(\omega)=\cases{\sum_{k=1}^\infty \sum_{C' \in \cD(N;k), \; C' \subset C} Z_{C'}(\omega) \b1_{C'}, &if $\omega\in \Omega^*$, \cr
                     0,                                  &if $\omega\notin \Omega^*$.\cr}$$

  The random variables $X^C$, $C \in \cD(N;1)$, are built from disjoint
collections of weights and none of these collections
 contains the weight $Z_{D^N}$. 
Thus the random variables $X^C$, $C \in \cD(N;1)$, are independent and jointly independent of $Z_{D^N}$.  Moreover, if $s \in C$ for some  $C \in \cD(N;1)$, then the $C(\rho D^N, K^d)$-valued random variable  
$X^C(s+\cdot)$ has the same law as $X^{\rho D^N}$, which in turn has the same distribution as $\rho X(\rho^{-1} \cdot)$.  In particular, if $|z| <1$,
 then the function $z \b1_C$ belongs to the support of $X^C$ and 
$z\b1_D^N+X^C$ has the same law as $X^C$. Thus, for $|z| < 1$ we have 
$$\eqalign{\bP\{ 0 \in z + X^C(C)\} &= \bP\{ 0 \in X^C(C)\} \cr
                                 &= \bP\{0 \in X^{\rho D^N}(\rho D^N)\} \cr
                                 &= \bP\{0 \in X(D^N)\} \cr
                                 &= h. \cr}$$
Also, if $|z|=1$, then from the isosceles triangle property we see that
$|z + X^C(t)| =1$ for all $t \in C$, because $|X^C(t)| \le q^{-1}$ for all $t \in C$. Thus, for $|z|=1$ we have $\bP\{ 0 \in z + X^C(C)\} = 0$.

   Now $V_1$ is just the number of balls $C \in \cD(N;1)$ for which
$0 \in Z_{D^N} + X^C(C)$.  From what have observed in the previous
paragraph, conditional on $|Z_{D^N}| < 1$ the random variable $V_1$ 
has a ${\rm bin}(q^N,h)$ distribution, whereas $V_1 = 0$ 
when $|Z_{D^N}|=1$. 
Hence,
$$\eqalign{(1-h) &= \bP\{\cZ = \emptyset\} \cr
                 &= \bP\{|Z_{D^N}| = 1\} + \bP\{|Z_{D^N}| < 1\} (1-h)^{q^N} \cr
                 &= (1-q^{-d}) + q^{-d}(1-h)^{q^N}. \cr}$$
Moreover, if we let $Q$ denote the distribution of a ${\rm bin}(q^N,h)$ random variable conditioned on being non-zero then $Q$ is the distribution of $V_1$ conditional on $\cZ \ne \emptyset$.  In particular,
$$\eqalign{\bP[V_1 | \cZ \ne \emptyset] &= \sum_{i=1}^{q^N} iQ(i)\cr
                    &= {{q^N h} \over {1 - (1-h)^{q^N}}}\cr
                    &= {{q^N h} \over {1 - q^d[(1-h)-(1-q^{-d})]}}\cr
                    &= q^{N-d}.\cr}$$

   Continuing this line of argument shows that, conditional on $\cZ \ne \emptyset$, the  distribution of the sequence of random variables $\{V_n\}_{n=0}^\infty$ is that of a supercritical Galton-Watson branching process with offspring distribution $Q$.  Applying Theorem I.8.1 in [Harris, 1963] then gives that  $\kappa_n(D^n)=q^{(d-N)n} V_n$ converges almost surely and in $\cL^2$ as $n \rightarrow \infty$ to some random variable with expectation $h$.

 An almost identical argument shows that $\kappa_n(C)$ converges almost
surely and in $\cL^2$ for each  
 $\Gamma(N)$.
As in the proof of Theorem 6.3, we can use the fact that the algebra of sets
consisting of finite unions of balls from $\Gamma(N)$ is a weak convergence
determining class  to conclude that $\kappa_n$ converges  weakly
almost surely as $n \rightarrow \infty$ to a random measure $\kappa$.  
It is clear that $\kappa$ is square-integrable and adapted. Again as
in the proof of Theorem 6.3, we may choose the limit  $\kappa$ 
to be homogeneous.  By construction, $\kappa$ is supported on $\cZ$.  As $\bP[\kappa(D^N)] = h$, the theorem follows from Lemma 7.1. \hfil\break\rightline{\endproof}

\vskip 24pt

\proclaim Corollary 7.3.  The closed support of $L_0$ coincides with $\cZ$ almost surely.

\vskip 18pt
\nobreak\noindent
{\it Proof.}  In the notation of Theorem 7.2, we see from Corollary I.10.3 of [ Athreya and Ney, 1972] that the conditional probability $\bP\{L_0(D^N) = 0 | \cZ \ne \emptyset\}$ is just the extinction probability for a branching process with offspring $Q$, and this latter probability is of course $0$
because $Q$ is supported on $\{1, \ldots, q^N\}$.  Similarly, for any other coset $C \in \Gamma(N)$ we see that $L_0(C) > 0$ almost surely whenever $\cZ \cap C \ne \emptyset$, and the result follows.  \hfil\break\rightline{\endproof}

\vskip 24pt

   The proof of the following Hausdorff measure intrinsic 
construction of $L_0$ may be found in [Evans, 1993]. The
proof involves precise estimates on the
tail  of the distribution of $L_0(C)$, $C \in \Gamma(N)$,
and, roughly speaking, the exact asymptotic rate of decrease 
as $n \rightarrow \infty$ for 
$L_0(t + \rho^n D^N)$ when $t$ is a ``typical'' point in $\cZ$.

\vskip 24pt

\noindent
{\bf Notation.}  Define a function $f:\{q^{-n}\}_{n=0}^\infty \rightarrow [0,\infty[$ by setting $f(r) = r^{N-d} (\log|\log r|)^{d/N}$.
In the notation of \S 6, let $f-m(\cdot)$ denote the Hausdorff measure on $D^N$ constructed using the measure function $f$.   

\vskip 24pt 

\proclaim Theorem 7.4.  For some constant $0<c<\infty$ we have $f-m(\cdot \cap \cZ) = cL_0$.

   The proof of Theorem 7.5 below, our analogue of Trotter's theorem
on the continuity of Brownian local time, may also be found in
[Evans, 1993] and we will just outline the main ideas.  As in
 the proof of Theorem 6.1,
$L_x$ is approximated by the random measure $B \rightarrow L_x^{(n)}(B) = \int_B \phi_n(X_t - x) \, dt$.  Because of the
martingale properties of this sequence, it is possible to
establish uniform bounds on the exponential moments of the random variables
$|L_x^{(n)}(B) - L_y^{(n)}(B)|$.
These bounds can then be fed into the standard machinery of the
Garsia - Rodemich - Rumsey lemma (see [Garsia et al., 1970] and
the generalisation in [Preston, 1971]) to obtain equicontinuity of the sequence
$\{L_\cdot^{(n)}(B)\}_{n=0}^\infty$.

\vskip 24pt
\proclaim Theorem 7.5.  There is a weakly continuous version of the 
measure-valued process $x \mapsto L_x$.  This version may be chosen so that $L_x$ is homogeneous for all $x \in D^d$.

\vskip 24pt
\proclaim Corollary 7.6. Suppose that $\mu$ is a finite measure on $D^d$.  Then  the unique integrable, adapted, homogeneous random measure  with characteristic measure $\mu$ may be represented as $\kappa(B) = \int L_x(B) \, \mu(dx)$, where $x \mapsto L_x$ is the version of the local time given by Theorem 7.5.

\vskip 18pt
\nobreak\noindent
{\it Proof.}  Theorem 7.5 guarantees that all the integrals $\int L_x(B) \, \mu(dx)$ are well-defined and that the resulting random measure is homogeneous.  By Fubini's theorem
$$\bP[\int \int f(X_t)\, L_x(dt) \, \mu(dx)] = \int f(x) \, \mu(dx)$$
for all nonnegative measurable functions $f$, and the result follows from Theorem 6.1.  \hfil\break\rightline{\endproof}      

\vskip 24pt
\proclaim Corollary 7.7.  The random set $X(D^N)$ has nonempty interior almost surely.

\vskip 18pt
\nobreak\noindent
{\it Proof.}  Note that $X(D^N)$ contains the open set $\{x : L_x(D^N) > 0\}$.  By Corollary 7.6 $\int_{D^d} L_x(D^N) \, dx = \int_{D^N} 1 \, dt$ almost surely and so  $\{x : L_x(D^N) > 0\} \ne \emptyset$ almost surely.  
\hfil\break\rightline{\endproof}

\vskip 24pt
\noindent
{\bf Remark.}  It is a consequence of the Ray-Knight theorem that
Brownian local time is almost surely strictly positive on the
interior of the range of linear Brownian motion.  We remark
without proof that, using techniques similar to those found
in the proof of Theorem 7.4 given in [Evans, 1993],
it is possible to show that
the interior of $X(D^N)$ coincides with
$\{x : L_x(D^N) > 0\}$ almost surely.

\vfill\eject
\centerline{\bf 8. Other $K$-Gaussian random series}

\vskip 24pt

  Recall the specific instance of the $(N,1)$ $K$-Brownian
motion that we built in \S 5.  Our construction was
a $C(D^N, K)$-valued random of the form $\sum_{n=0}^\infty Y_n f_n$,
where $\{Y_n\}_{n=0}^\infty$ is a sequence of independent $K$-valued,
$K$-Gaussian random variables and $\{f_n\}_{n=0}^\infty \subset C(D^N, K)$.
In this section we will consider some other examples of such random
series. 

Unlike the Euclidean case, where necessary and sufficient conditions for the almost sure convergence of Gaussian random series are often rather delicate (cf. [Marcus and Pisier, 1981]), the corresponding question in our setting  
is almost trivial.

\vskip 24pt
\noindent
\proclaim Lemma 8.1.
Suppose that $(E, \Vert~~~\Vert_E )$ is a Banach space.
Consider $\{ f_n \}_{n=0}^{\infty} \subset E$ and a sequence
$\{ Y_n \}_{n=0}^{\infty}$ of independent, $K$-valued, $K$-Gaussian random variables.
The series $\sum_{n=0}^{\infty} Y_n f_n$ converges almost surely in $E$ if and
only if $\Vert Y_n \Vert_{\infty} \Vert f_n \Vert_E \rightarrow 0$ as $n \rightarrow \infty$.  If the series converges, then the limit is an
$E$-valued, $K$-Gaussian random variable.

\vskip 18pt
\noindent
{\it Proof.}
From the ultrametric inequality we have that the series converges if and
only if $\vert Y_n \vert \Vert f_n \Vert_E \rightarrow 0$ almost surely as $n \rightarrow \infty$.
From the Borel-Cantelli lemmas we see that $\vert Y_n \vert \Vert f_n \Vert_E \rightarrow 0$
if and only if
$$
\sum_{n=0}^{\infty} \bP\{ \vert Y_n \vert \Vert f_n \Vert_E > \epsilon \} < \infty \leqno(8.1.1)
$$
for all $\epsilon >0$.
Since $\bP\{ \vert Y_n \vert \Vert f_n \Vert_E > \epsilon \}=0$ when
$\epsilon \ge \Vert Y_n \Vert_{\infty} \Vert f_n \Vert_E$ and
$\bP\{ \vert Y_n \vert \Vert f_n \Vert_E > \epsilon \} \ge 1-q^{-1}$ when
$\epsilon < \Vert Y_n \Vert_{\infty} \Vert f_n \Vert_E$, it is clear
that (8.1.1) occurs if and only if $\Vert Y_n \Vert_{\infty} \Vert f_n \Vert_E > \epsilon$
for finitely many $n$, and so the result on convergence follows.

   It follows from part (ii) of Theorem 4.8 and Theorem 4.9 that
$\sum_{n=0}^m Y_n f_n$ is an $E$-valued, $K$-Gaussian random variable.
Corollary 4.5 shows that $\sum_{n=0}^\infty Y_n f_n$ is also an $E$-valued, $K$-Gaussian random variable when the series converges.
\hfil\break\rightline{\endproof}

\vskip 24pt
For the remainder of this section, we will be concerned with random series
in the case when $K= \bQ_p$ and our Banach space is $C( \bZ_p, \bQ_p )$
equipped with the supremum norm $\| \; \|_C$.
In particular we will investigate the problem of representing stationary
$C( \bZ_p, \bQ_p )$-valued,
$\bQ_p$-Gaussian random variables
 as random series. (As in  \S 5, we say that a 
$C( \bZ_p, \bQ_p )$-valued,
$\bQ_p$-Gaussian random variable $X$ is stationary if $X(\cdot + s)$ has the same law
as $X$ for all $s \in \bZ_p$.)
Our first result concerns the Mahler basis introduced in Example 3.9.

\vskip 24pt
\noindent
{\bf Note.} {\it For this section the notation $X$ will no longer be reserved
exclusively for $K$-Brownian motion.}

\vskip 24pt
\noindent
{\bf Definition.}
Let $\{ Z_n \}_{n=0}^{\infty}$ be a sequence of independent, $\bQ_p$-valued,
$\bQ_p$-Gaussian random variables such that $\Vert Z_n \Vert_{\infty} =1$
for all $n$.
Suppose that $\{ a_n \}_{n=0}^{\infty} \subset \bQ_p$ is such that
$\vert a_n \vert \rightarrow 0$ as $n \rightarrow \infty$.
We say that the $C( \bZ_p, \bQ_p )$-valued, $\bQ_p$-Gaussian random variable
 $X$ defined by $X(t)= \sum_{n=0}^{\infty} a_n Z_n {t \choose n}$,
$t \in \bZ_p$, is a random Mahler series.

\vskip 24pt
Since $\Vert {\cdot \choose n} \Vert_C =1$, we have from
Lemma 8.1 that $X$ is a well-defined $C( \bZ_p, \bQ_p )$-valued,
$\bQ_p$-Gaussian random variable.

\vskip 24pt
\noindent
\proclaim Theorem 8.2.
If $X= \{ \sum_{n=0}^{\infty} a_n Z_n {\cdot \choose n} \}$ is a random Mahler
series, then $X$ is stationary  if and only if $\vert a_n \vert \ge \vert a_{n+1} \vert$ for all $n$.

\vskip 18pt
\noindent
{\it Proof.}
Since $\bN$ is dense in $\bZ_p$, we have for any $s \in \bZ_p$
that there exists a sequence $\{m_k\}_{k=0}^\infty \in \bN$ such that 
$s = \lim_{k \rightarrow \infty} m_k$, and hence 
$\lim_{k \rightarrow \infty} \|X(\cdot + s) - X(\cdot + m_k)\|_C = 0$
almost surely.  Thus $X$ will be stationary if and only if $X(\cdot + m)$
has the same law as $X$ for all $m \in \bN$ which will, in turn occur
if and only if $Y = X(\cdot + 1)$ has
the same law as $X$.

From the calculation on p. 152 of [Schikhof, 1984] (essentially
a consequence of
the ``Pascal's triangle'' recurrence relation for the binomial
coefficients), we find that
$$
Y(t) = \sum_{n=0}^{\infty} (a_n Z_n + a_{n+1} Z_{n+1} ) {t \choose n} .
$$
It is clear from Theorem 4.9 that the finite dimensional distributions
of the sequence $\{ a_n Z_n + a_{n+1} Z_{n+1} \}_{n=0}^{\infty}$ are $\bQ_p$-Gaussian.
Hence $X$ will be stationary if and only if
$$
\Vert a_n Z_n + a_{n+1} Z_{n+1} \Vert_{\infty} = \Vert a_n Z_n \Vert_{\infty} \leqno(8.2.1)
$$
for all $n$, and the sequence $\{ a_n Z_n + a_{n+1} Z_{n+1} \}_{n=0}^{\infty}$ is independent.
Observe from part (i) of Theorem 4.8 that, as $Z_n$ and $Z_{n+1}$ are independent
they are orthogonal in $\cL^\infty(\bP)$, and so (8.2.1) is equivalent to requiring that
$$
\vert a_n \vert \vee \vert a_{n+1} \vert = \vert a_n \vert . \leqno(8.2.2)
$$

Suppose that $X$ is stationary, then (8.2.2) implies that
$\vert a_n \vert \ge \vert a_{n+1} \vert$ for all $n$.

Conversely, suppose that $\vert a_n \vert \ge \vert a_{n+1} \vert$ for all $n$.
Then (8.2.2) holds. Applying Theorem 4.11, we see that  the sequence $\{ a_n Z_n + a_{n+1} Z_{n+1} \}_{n=0}^{\infty}$ is independent if  for each $n \in \bN$ the collection of vectors
$\alpha_0 =(a_0 , a_1 , 0, \ldots ,0)$,
$\alpha_1 =(0,a_1 , a_2 , 0, \ldots ), \ldots , \alpha_n =(0, \ldots , 0, a_n , a_{n+1} )$
is orthogonal in $( \bQ_p )^{n+1}$.
However, for $\lambda_0 , \ldots , \lambda_n \in \bQ_p$, an induction based
on Lemma A.1 in the Appendix shows that
$$
\eqalign{\vert \lambda_0 \alpha_0&+ \ldots + \lambda_n \alpha_n \vert \cr
&= \vert \lambda_0 \vert \vert a_0 \vert \vee \vert \lambda_0 + \lambda_1 \vert
\vert a_1 \vert \vee \ldots \vee \vert \lambda_{n-1} + \lambda_n \vert \vert a_n \vert  \vee \vert \lambda_n \vert \vert a_{n+1} \vert \cr
&= \vert \lambda_0 \vert \vert a_0 \vert \vee \ldots \vee \vert \lambda_n \vert \vert a_n \vert \cr
&= \vert \lambda_0 \vert \vert \alpha_0 \vert \vee \ldots \vee \vert \lambda_n \vert \vert \alpha_n \vert , \cr}
$$
and so $\alpha_0 , \ldots , \alpha_n$ are orthogonal, as required.
\hfil\break\rightline{\endproof}

\vskip 24pt
A remarkable feature of the Gaussian theory is that stationary processes
on the circle can be represented as random Fourier series with independent
Fourier coefficients.
With this in mind one might hope that all 
$C( \bZ_p, \bQ_p )$-valued, $\bQ_p$-Gaussian
random variables have the form given in Theorem 8.2.
In Corollary 8.4 below, we show that not only is this not the case, but in
fact there is no orthonormal basis for $C( \bZ_p, \bQ_p )$ that ``works.''
First, however, we obtain a result similar to Theorem 8.2 for the van der Put basis
$\{ e_n \}_{n=0}^{\infty}$ introduced in Example 3.10.

\vskip 24pt
\noindent
{\bf Definition.}
Let $\{ Z_n \}_{n=0}^{\infty}$ be a sequence of independent
$\bQ_p$-valued, $\bQ_p$-Gaussian random variables such that
$\Vert Z_n \Vert_{\infty} =1$ for all $n$.
Suppose that $\{ a_n \}_{n=0}^{\infty} \subset \bQ_p$ is such that
$\vert a_n \vert \rightarrow 0$ as $n \rightarrow \infty$.
We say that the $C( \bZ_p, \bQ_p )$-valued, $\bQ_p$--Gaussian
random variable $X$ defined by
$X = \sum_{n=0}^{\infty} a_n Z_n e_n$
is a random van der Put series.

As in the remarks following the definition of Mahler series, we see that $X$ is 
a well-defined $C(\bZ_p, \bQ_p)$-valued, $\bQ_p$-Gaussian random variable.

\vskip 24pt
\noindent
\proclaim Theorem 8.3.
If $X= \{ \sum_{n=0}^{\infty} a_n Z_n e_n\}$ is a random van der
Put series, then $X$ is stationary if and only if
$$
\vert a_0 \vert \ge \vert a_1 \vert \ge \vert a_p \vert \ge \ldots \ge
\vert a_{p^n} \vert \ge \vert a_{p^{n+1}} \vert \ge \ldots \leqno(8.3.1)
$$
and
$$
\vert a_{p^n} \vert = \vert a_{p^n +1} \vert = \ldots = \vert
a_{p^{n+1} -1} \vert \leqno(8.3.2)
$$
for all $n$.

\vskip 18pt
\noindent
{\it Proof.}
As in the proof of Theorem 8.2, we have that if we set
$Y(t)=X(t+1)$, $t \in \bZ_p$, then $X$ will be stationary if and
only if $Y$ has the same law as $X$.

From Exercise 62.F in [Schikhof, 1984], we find that
$$
Y(t)= \sum_{n=0}^{\infty} B_n e_n (t)
$$
where
$$
B_n = \cases{a_0 Z_0 + a_1 Z_1&if $n=0$,\cr
a_{n+1} Z_{n+1} -a_{p^s} Z_{p^s}&if $n=rp^s -1 ,~ s \in \bN , ~ 2 \le r \le p$,\cr
a_{n+1} Z_{n+1}&otherwise.\cr}
$$
It is clear from Theorem 4.9 that the finite dimensional distributions
of the sequence $\{ B_n \}_{n=0}^{\infty}$ are $\bQ_p$-Gaussian and so $X$ will be stationary if and only if
$$
\Vert B_n \Vert_{\infty} = \vert a_n \vert \leqno(8.3.3)
$$
for all $n$ and the sequence $\{ B_n \}_{n=0}^{\infty}$ is independent.

Suppose that the conditions (8.3.1) and (8.3.2) holds.
We have from part (i) of Theorem 4.8 that the sequence
$\{Z_n\}_{n=0}^\infty$ is orthonormal in $\cL^\infty(\bP)$, and hence
 (8.3.3) holds.
Since $\{ B_0 , \ldots , B_{p-1} \}$ is contained
in the linear span of $\{ Z_0 , \ldots , Z_p \}$ in $\cL^\infty(\bP)$ and
and $\{ B_{p^s} , \ldots , B_{p^{s+1} -1} \}$
is contained in the linear span of
$\{ Z_{p^s +1} , \ldots , Z_{p^{s+1}} \}$  in $\cL^\infty(\bP)$ for
$s=1,2, \ldots$, it 
suffices to show that the elements within each of these subsets of
$\{B_n\}_{n=0}^\infty$ are independent.  This is in turn equivalent, again by part (i) of Theorem 4.8, to showing that each
such subset is orthogonal.
From the orthonormality of $\{Z_n\}_{n=0}^\infty$ and Lemma A.2 in the Appendix, we have
$$
\eqalign{\Vert \lambda_0 B_0&+ \lambda_1 B_1 + \ldots + \lambda_{p-1} B_{p-1} \Vert_{\infty} \cr
&= \Vert \lambda_0 (a_0 Z_0 + a_1 Z_1 ) + \lambda_1 (a_2 Z_2 -a_1 Z_1 )+ \ldots
+ \lambda_{p-1} (a_p Z_p -a_1 Z_1 ) \Vert_{\infty} \cr
&= \vert \lambda_0 \vert \vert a_0 \vert \vee ( \vert \lambda_0 -( \lambda_1 + \ldots + \lambda_{p-1} ) \vert  \vee \vert \lambda_1 \vert \vee \ldots \vee \vert \lambda_{p-1} \vert )
\vert a_1 \vert \vee \vert \lambda_{p-1} \vert \vert a_p \vert \cr
&= \vert \lambda_0 \vert \vert a_0 \vert \vee ( \vert \lambda_1 \vert \vee \ldots \vee \vert
\lambda_{p-1} \vert ) \vert a_1 \vert \cr
&= \vert \lambda_0 \vert \Vert B_0 \Vert_{\infty} \vee \vert \lambda_1 \vert
\Vert B_1 \Vert_{\infty} \vee \ldots \vee \vert \lambda_{p-1} \vert \Vert B_{p-1} \Vert_{\infty} ,\cr}
$$
so that $B_0 , \ldots B_{p-1}$ are orthogonal.
A similar argument using Lemma A.3 in the Appendix establishes that $B_{p^s} , \ldots , B_{p^{s+1} -1}$
are orthogonal for each $s=1,2, \ldots$ and hence completes the proof that $X$ is stationary.

Conversely, if $X$ is stationary, then by  the orthonormality of $\{Z_n\}_{n=0}^\infty$ the condition 8.3.3
is equivalent to requiring that
$$
\eqalign{\vert a_0 \vert&= \vert a_0 \vert \vee \vert a_1 \vert ,\cr
\vert a_n \vert&= \vert a_{n+1} \vert \vee \vert a_{p^s} \vert~{\rm if}~
n=rp^{-s} ,~ s \in \bN , ~ 2 \le r \le p ,\cr}
$$
and
$$
\vert a_n \vert = \vert a_{n+1} \vert~{\rm for~all~other}~n .
$$
It is straightforward to check that this implies (8.3.1) and (8.3.2).
\hfil\break\rightline{\endproof}

\vskip 24pt
 The 
statement of the following result is
similar to the statement of
Corollary 9.7 of [Evans, 1989a], with
the exception that for the former result it is
only assumed that the sequence $\{ f_n \}_{n=0}^{\infty}$
is linearly independent.
We thank John Taylor for pointing
out that the proof in [Evans, 1989a]
implicitly uses something stronger than this
purely algebraic condition.  
We do not know if another proof can be given or whether a counterexample exist.

\vskip 24pt
\noindent
\proclaim Corollary 8.4.
There is no orthonormal basis $\{ f_n \}_{n=0}^{\infty}$ for
$C( \bZ_p, \bQ_p )$
such that every stationary 
$C(\bZ_p,\bQ_p)$-valued, $\bQ_p$-Gaussian random  variable $X$
is of the form $X = \sum_{n=0}^{\infty} A_n f_n$ for some sequence
$\{ A_n \}_{n=0}^{\infty}$ of independent, $\bQ_p$-valued, $\bQ_p$-Gaussian
random variables.

\vskip 18pt
\noindent
{\it Proof.}
Suppose that $\{ f_n \}_{n=0}^{\infty}$ has the requisite properties.

For $m\in \bN$, set
$X = \sum_{n=0}^{p^m -1} B_n e_n$,
where $\{ e_n \}_{n=0}^{\infty}$ is the van der Put basis and
$\{ B_k \}_{k=0}^{p^m -1}$ is a set of independent $\bQ_p$-Gaussian
random variables with
$$
\Vert B_0 \Vert_{\infty} = \Vert B_1 \Vert_{\infty} = \ldots = \Vert
B_{p^n -1} \Vert_{\infty} =1 .
$$
From Theorem 8.3, we see that $X$ is stationary.
By assumption we have that
$X = \sum_{n=0}^{\infty} A_n f_n$
where $\{ A_n \}_{n=0}^{\infty}$ is a set of independent, $\bQ_p$-Gaussian
random variables.

   As the sequence $\{ e_n \}_{n=0}^{\infty}$ (resp. $\{ f_n \}_{n=0}^{\infty}$)
is linearly independent, each random variable $B_n$ (resp. $A_n$)
is obtained by the application of some continuous linear functional to $X$.
Thus the closed linear span in $\cL^\infty(\bP)$ of $\{B_n\}_{n=0}^{p^m-1}$
(resp. $\{A_n\}_{n=0}^\infty$) is contained in the closed linear span
of the set $\{T(X) : T \in C( \bZ_p, \bQ_p )^*\}$.  The two
reverse containments  obviously hold, and so the closed linear span in $\cL^\infty(\bP)$ of $\{B_n\}_{n=0}^{p^m-1}$ and the closed linear span in $\cL^\infty(\bP)$ of $\{A_n\}_{n=0}^\infty$ are equal.

By part (i) of Theorem 4.8, the sequence $\{ B_0 , \ldots , B_{p^m -1} \}$
(resp. $\{ A_n \}_{n=0}^{\infty}$) is orthogonal.  In particular,
$\{ B_0 , \ldots , B_{p^m -1} \}$ (resp. 
$\{ A_n : A_n \ne 0 \}$) is linearly independent.  As
the closed linear spans in $\cL^\infty(\bP)$ of $\{B_n\}_{n=0}^{p^m-1}$
and  $\{ A_n : A_n \ne 0 \}$ coincide, we find by equating dimensions
that the set $\{ A_n : A_n \ne 0 \}$
must have exactly $p^m$ elements.  We will index these as $A_{m,n}$
for $n=0, \ldots,p^m-1$.  Then
$X = \sum_{n=0}^{p^m -1} A_{m,n}f_{m,n},$
where $\{ f_{m,0},  \ldots , f_{m,p^m -1}\} \subset \{ f_n \}_{n=0}^{\infty}$.

Using the independence of $B_0 , \ldots , B_{p^m -1}$
we can write the closed support of the law of the $C( \bZ_p, \bQ_p )$-valued
random variable $X$ as
$$
\{ \sum_{n=0}^{p^m -1} b_n e_n : \vert b_0 \vert \le  1, \ldots , \vert
b_{p^m -1} \vert \le 1 \}.
$$
Using the independence of $A_{m,0} , \ldots , A_{m,p^m -1}$ this closed
support is also

$$
\{ \sum_{n=0}^{p^m -1} a_n f_{m,n} : \vert a_0 \vert \le \Vert A_{m,0} \Vert , \ldots ,
\vert a_{p^m -1} \vert \le \Vert A_{m,p^m -1} \Vert \} .
$$
Any element of the linear span of $\{e_0, \ldots, e_{p^m-1}\}$ can
be written in the form $d\sum_{n=0}^{p^m -1} b_n e_n$ where $d \in \bQ_p$,
$\vert b_0 \vert \le  1, \ldots , \vert
b_{p^m -1} \vert \le 1$.
Any element of the linear span of $\{f_{m,0}, \ldots, f_{m,p^m-1}\}$
can be written as
$c\sum_{n=0}^{p^m -1} a_n f_{m,n}$ where $c \in \bQ_p$,
$\vert a_0 \vert \le \Vert A_{m,0} \Vert , \ldots ,
\vert a_{p^m -1} \vert \le \Vert A_{m,p^m -1} \Vert $.
Thus these two linear spans are equal.  
Consequently,
the linear span of $\{ e_n : n=0,1, \ldots \}$ is
the same as the linear span of 
$\bigcup_{m=0}^{\infty} \{f_{m,n} : n=0 , \ldots
, p^m -1 \}$.  

There can be no basis function $f_N$
not in the set $\bigcup_{m=0}^{\infty} \{f_{m,n} : n=0 , \ldots
, p^m -1 \}$, because such a function would (by orthonormality)
be at distance $1$ from every function in the closed linear span
of the set,  and the closed linear span of $\{ e_n : n=0,1, \ldots \}$
is all of $C(\bZ_p, \bQ_p)$.  Thus every basis function $f_N$
belongs to the linear span of  $\{ e_n : n=0,1, \ldots \}$.
In particular, every basis function $f_N$ is locally constant.

We can essentially repeat the argument we have gone through up to now
with the van der Put basis replaced by the Mahler basis and Theorem 8.3
replaced by Theorem 8.2. This allows us to conclude that every basis function
$f_N$ is in the linear span of the Mahler functions, and hence is
a polynomial.

The only locally constant polynomials are the constant functions. Any
orthonormal basis of $C(\bZ_p, \bQ_p)$ contains at most one
constant function, so we obtain a contradiction.
\hfil\break\rightline{\endproof}

\vskip 24pt
\noindent
{\bf Remark.} Given Corollary 8.4, it is natural to ask if there
is a simple characterisation of the stationary
$C(\bZ_p, \bQ_p)$-valued, $\bQ_p$-Gaussian random variables that
can be represented as Mahler or van der Put series.  In particular,
one can ask if the $(1,1)$ $\bQ_p$-Brownian motion is a Mahler
or van der Put series.  We don't know.  Also, one can ask
if there is any concrete method for constructing all the stationary
$C(\bZ_p, \bQ_p)$-valued, $\bQ_p$-Gaussian random variables.
In a forthcoming paper, we discuss the representation of
stationary random variables as ``moving average stochastic integrals''
against ``$\bQ_p$-white noise''.  We show that $\bQ_p$-Brownian
motion may be represented in this manner.  However, we don't
know if all stationary $C(\bZ_p, \bQ_p)$-valued, $\bQ_p$-Gaussian random variables have such a representation,
although we expect that the answer is negative.

\vfill\eject

\vskip 30pt

\centerline{\bf Appendix}

\vskip 24pt

   The following three lemmas were used in \S 8.

\vskip 24pt
\noindent
\proclaim Lemma A.1.
Suppose that $a \ge b \ge 0$.
If $\alpha , \beta \in K$, then
$$
( \vert \alpha \vert a) \vee ( \vert \alpha + \beta \vert b)=
( \vert \alpha \vert a ) \vee ( \vert \beta \vert b) .
$$

\vskip 18pt
\noindent
{\it Proof.}
Suppose, first of all, that $\vert \alpha \vert a \ge \vert \beta \vert b$.
From the ultrametric inequality, we have that
$$
\vert \alpha + \beta \vert b \le ( \vert \alpha \vert \vee \vert
\beta \vert )b \le \vert \alpha \vert a ,
$$
and the equality holds.

On the other hand, if $\vert \alpha \vert a < \vert \beta \vert b$,
then $\vert \alpha \vert < \vert \beta \vert$.
Applying the isosceles triangle property, we have that $\vert \alpha + \beta \vert = \vert \beta \vert$,
and the equality also holds.
\hfil\break\rightline{\endproof}

\vskip 24pt
\noindent
\proclaim Lemma A.2.
Suppose that $a \ge b \ge c \ge 0$.
If $\alpha_1 , \ldots , \alpha_n \in K$, then
$$
\vert \alpha_1 \vert a \vee ( \vert \alpha_1 - ( \alpha_2 + \ldots + \alpha_n ) \vert
\vee \vert \alpha_2 \vert \vee \ldots \vee \vert \alpha_{n-1} \vert ) b \vee
\vert \alpha_n \vert c
$$
$$
= \vert \alpha_1 \vert a \vee ( \vert \alpha_2 \vert \vee \ldots \vee
\vert \alpha_n \vert )b .
$$

\vskip 18pt
\noindent
{\it Proof.}
The left-hand side is at most the right-hand side by the ultrametric inequality.
The reverse inequality is clear except in the following cases.
For ease of notation, we set $\beta = \alpha_1 - ( \alpha_2 + \ldots + \alpha_n )$.

{\bf Case I.}
$$
\vert \alpha_1 \vert a \vee ( \vert \alpha_2 \vert \vee \ldots \vee \vert
\alpha_{n-1} \vert )b < \vert \beta \vert b \leqno(A.2.1)
$$
$$
\vert \alpha_n \vert c \le \vert \beta \vert b \leqno(A.2.2)
$$

{\bf Case II.}
$$
\vert \alpha_1 \vert a \vee ( \vert \beta \vert \vee \vert \alpha_2 \vert \vee
\ldots \vee \vert \alpha_{n-1} \vert )b< \vert \alpha_n \vert c . \leqno(A.2.3)
$$

Suppose that Case I holds.
From (A.2.1) and the ultrametric inequality, we have that
$$
\eqalign{\vert \alpha_1 -( \alpha_2 + \ldots + \alpha_{n-1} ) \vert & \le ( \vert \alpha_1
\vert \vee \ldots \vee \vert \alpha_{n-1} \vert ) \cr
& < \vert \beta \vert \cr
& = \vert ( \alpha_1 -( \alpha_2 + \ldots + \alpha_{n-1} )) - \alpha_n \vert .\cr}
$$
The isosceles triangle property then implies that
$$
\vert \alpha_n \vert < \vert \beta \vert ,
$$
so the reverse inequality holds.

Since
$$
\vert \alpha_n \vert \le \vert \alpha_1 \vert \vee \ldots \vee \vert
\alpha_{n-1} \vert \vee \vert \beta \vert
$$
by the ultrametric inequality, we see that Case II cannot hold.
\hfil\break\rightline{\endproof}

\vskip 24pt
\noindent
\proclaim Lemma A.3.
Suppose that $a \ge b \ge 0$.
If $\alpha_1 , \ldots , \alpha_n \in K$, then
$$
( \vert \alpha_1 \vert \vee \ldots  \vert \alpha_{n-1} \vert )
a \vee \vert \alpha_n \vert b \vee \vert \alpha_1 + \ldots + \alpha_n \vert a
$$
$$
=( \vert \alpha_1 \vert \vee \ldots \vee \vert \alpha_n \vert )a .
$$

\vskip 18pt
\noindent
{\it Proof.}
The left-hand side is at most the right-hand side by the ultrametric inequality,
and the reverse inequality is clear except when $\vert \alpha_n \vert > ( \vert \alpha_1 \vert \vee \ldots \vee \vert \alpha_{n-1} \vert )$.
In this case, we have from the ultrametric inequality that
$\vert \alpha_n \vert > \vert \alpha_1 + \ldots + \alpha_{n-1} \vert$,
and so, by the isosceles triangle property,
$$
\vert \alpha_1 + \ldots + \alpha_n \vert = \vert \alpha_n \vert ,
$$
and the result follows.
\hfil\break\rightline{\endproof}

\vskip 24pt

 The following result was used in the proof of Theorem 12.4 and follows from a straightforward monotone class argument.

\vskip 24pt

\proclaim Lemma A.4.  Let $(\Sigma,{\cal A},{\bf Q})$ be  a probability space.  Suppose that we have ${\cal A}={\cal A}_1\vee{\cal A}_2\vee{\cal A}_3\vee{\cal A}_4$, where
${\cal A}_1, {\cal A}_2, {\cal A}_3, {\cal A}_4$ are independent sub-$\sigma$-fields. For a ${\bf Q}$-integrable random variable $Y$ we have
$${\bf Q}(Y|{\cal A}_1)={\bf Q}({\bf Q}(Y|{\cal A}_1\vee{\cal A}_2)|{\cal A}_1\vee{\cal A}_3).$$

\vskip 24pt

   The next result, which was used also used in the proof of Theorem 12.4, is just a discrete time version of Meyer's energy inequality.

\vskip 24 pt

\proclaim Lemma A.5.  Let $(\Sigma,{\cal A},{\bf Q})$ be  a probability space.  Suppose that ${\cal A}_1\subset{\cal A}_2\subset\cdots\subset{\cal A}_n$ are sub-$\sigma$-fields of ${\cal A}$
and $Y_1, Y_2,\ldots, Y_n$ are non-negative, square-integrable random variables. Then
$${\bf Q}([\sum_k{\bf Q}(Y_k|{\cal A}_k)]^2)\le4{\bf Q}([\sum_k Y_k]^2).$$

\vskip 18pt \noindent
{\it Proof.}   We have
$$\eqalign{{\bf Q}([\sum_k{\bf Q}(Y_k|{\cal A}_k)]^2)
   &\le2\sum_k\sum_{l\ge k}{\bf Q}({\bf Q}(Y_k|{\cal A}_k){\bf Q}(Y_l|{\cal A}_l))\cr
   &=2\sum_k\sum_{l\ge k}{\bf Q}({\bf Q}(Y_k|{\cal A}_k)Y_l)\cr
   &\le2{\bf Q}(\sum_k{\bf Q}(Y_k|{\cal A}_k)\sum_l Y_l)\cr
   &\le2{\bf Q}([\sum_k{\bf Q}(Y_k|{\cal A}_k)]^2)^{1\over 2}{\bf Q}([\sum_l Y_l]^2)^{1\over 2},\cr}$$
where the last line follows from the Cauchy-Schwarz inequality.
\hfil\break\rightline{\endproof}

\vfill\eject

\vskip 30pt

\centerline{\bf References}

\vskip 24pt\noindent
Adler, R. (1990). {\it An Introduction to Continuity, Extrema, and Related Topics for General Gaussian Processes.} Institute of Mathematical Statistics.
\bigskip\noindent
Albeverio, S. and Karwowski, W. (1991). Diffusion on $p$-adic numbers. In {\it Gaussian Random Fields (Nagoya 1990).} World Scientific.
\bigskip\noindent
Albeverio, S. and Karwowski, W. (1994). A random walk on $p$-adics - the
generator and its spectrum. {\it Stochastic Process. Appl.} {\bf 53} 1-22.
\bigskip\noindent
Athreya, K.B. and Ney, P.E. (1972). {\it Branching Processes.} Springer.
\bigskip\noindent
Barlow, M.T., Perkins, E.A. and Taylor, S.J. (1986). The behaviour and construction of local times for L\'evy processes.  In {\it Seminar on Stochastic Processes 1984} (E. Cinlar, K.L. Chung, R.K. Getoor eds.) Birkh\"auser.
\bigskip\noindent
Brillinger, D. (1991). Some asymptotics of finite Fourier transforms of a stationary p-adic process. {\it J. Combinatorics and System Sciences} {\bf 16} 155-169.
\bigskip\noindent
Brydges, D., Evans, S.N. and Imbrie, J. (1992). Self-avoiding walk on the hierarchical lattice in four dimensions. {\it Ann. Probab.} {\bf 20} 82-124.
\bigskip\noindent
Cassels, J.W.S. (1986). {\it Local Fields.} Cambridge University Press.
\bigskip\noindent
Cambanis, S. and Rajput B.S. (1973).
Some zero-one laws for Gaussian processes.
{\it Ann. Probab.} {\bf 1} 304-312.
\bigskip\noindent
Cuoco, A.A. (1991). Visualizing the $p$-adic integers. {\it Amer. Math. Monthly} {\bf 98} 355-364.
\bigskip\noindent
Curtis, C.W. and Reiner, I. (1962). {\it Representation Theory of Finite Groups and Associative Algebras.} Wiley.
\bigskip\noindent
Doob, J.L. (1984). {\it Classical Potential Theory and Its Probabilistic
Counterpart.} Springer.
\bigskip\noindent
Dudley, R.M. (1989). {\it Real Analysis and Probability.} Wadsworth.
\bigskip\noindent
Dwork, B.M. (1982). {\it Lectures on p-adic Differential Equations.} Springer.
\bigskip\noindent
Dynkin, E.B. (1981). Additive functionals of several time-reversible Markov processes. {\it J. Funct. Anal.} {\bf 42} 64-101.
\bigskip\noindent
Ethier, S.N. and Kurtz, T.G. (1986). {\it Markov Processes: Characterization
and Convergence.} Wiley.
\bigskip\noindent
Evans, S.N. (1988a). Continuity properties of  Gaussian stochastic processes indexed by a local field. {\it Proc. London Math. Soc.} {\bf 56} 380-416. 
\bigskip\noindent
Evans, S.N. (1988b). Sample path properties of Gaussian stochastic processes indexed by a local field. {\it Proc. London Math. Soc.} {\bf 56} 580-624.
\bigskip\noindent
Evans, S.N. (1989a). Local field Gaussian measures. In {\it Seminar on Stochastic Processes 1988} (E. Cinlar, K.L. Chung, R.K. Getoor eds.) Birkh\"auser.
\bigskip\noindent
Evans, S.N. (1989b). Local properties of L\'evy processes on a totally disconnected group. {\it J. Theoret. Probab.} {\bf 2} 209-259.
\bigskip\noindent
Evans, S.N. (1991). Equivalence and perpendicularity of local field Gaussian measures. In {\it Seminar on Stochastic Processes 1990} (E. Cinlar ed.) Birkh\"auser. 
\bigskip\noindent
Evans, S.N. (1992). Polar and non-polar sets for a tree indexed process. {\it Ann. Probab.} {\bf 20} 579-590.
\bigskip\noindent
Evans, S.N. (1993). Local field Brownian motion. {\it J. Theoret. Probab.}
{\bf 6} 817-850.
\bigskip\noindent
Feldman, J. (1958).
Equivalence and perpendicularity of Gaussian processes.
{\it Pacific J. Math.} {\bf 4} 699-708.
\bigskip\noindent
Fernique, X. (1975).
{\it Lecture Notes in Mathematics}, no. 480.
Springer.
\bigskip\noindent
Fitzsimmons, P.J. (1987). Homogeneous random measures and a weak order for the excessive measures of a Markov process. {\it Trans. Amer. Math. Soc.} {\bf 303} 421-478.
\bigskip\noindent
Fitzsimmons, P.J. and Salisbury, T.S. (1989). Capacity and energy for multiparameter Markov processes. {\it Ann. Inst. Henri Poincar\'e} {\bf 25} 325-350.
\bigskip\noindent
Garsia, A., Rodemich, E. and Rumsey Jr., H. (1970). A real variable lemma and the continuity of paths of some Gaussian processes. {\it Indiana Univ. Math. J.} {\bf 20} 565-578.
\bigskip\noindent
Geman, D. and Horowitz, J. (1980). Occupation densities. {\it Ann. Probab.} {\bf 8} 1-67.
\bigskip\noindent
Gregory, R.T. and Krishnamurthy, E.V. (1984). {\it Methods and Applications of Error-Free Computation.} Springer.
\bigskip\noindent
Guimier, F. (1989). Simplicit\'e du spectre de Liapounoff d'un produit de matrices al\'eatoires sur un corps ultram\'etrique. {\it C. R. Acad. Sci. Paris, Serie 1} {\bf 309} 885-889.
\bigskip\noindent
H\'ajek, J. (1959).
On a simple linear model in Gaussian processes.
In {\it Trans. Second Prague Conf. Information Theory}, 185-197.
\bigskip\noindent
Harris, T.E. (1963). {\it The Theory of Branching Processes.} Springer.
\bigskip\noindent
Iwasawa, K. (1972). {\it Lectures on p-adic L-functions.} Princeton University Press.
\bigskip\noindent
Jain, N.C. (1971).
A zero-one law for Gaussian processes.
{\it Proc. Amer. Math. Soc.} {\bf 29} 585-587.
\bigskip\noindent
Jain, N.C, and  Marcus, M.B. (1978).
Continuity of sub-gaussian processes, in {\it Advances in Probability}, Vol. 4. Marcel Dekker.
\bigskip\noindent
Kakutani, S. (1944a). On Brownian motion in $n$-spaces. {\it Proc. Imp. Acad. Tokyo} {\bf 20} 648-652.
\bigskip\noindent
Kakutani, S. (1944b). Two-dimensional Brownian motion and harmonic functions. {\it Proc. Imp. Acad. Tokyo} {\bf 20} 706-714.
\bigskip\noindent
Kallianpur, G. (1970).
Zero-one laws for Gaussian processes.
{\it Trans. Amer. Math. Soc.} {\bf 149} 199-211.
\bigskip\noindent
Khrennikov, A. Yu. (1990). Mathematical methods in non-Archmidean physics. {\it Russian Math. Surveys} {\bf 45} 87-125.
\bigskip\noindent
Kingman, J.F.C. (1973). An intrinsic description of local time. {\it J. London Math. Soc.} {\bf 6} 725-731.
\bigskip\noindent
Koblitz, N. (1980). {\it p-adic Analysis: a Short Course on Recent Work.} Cambridge University Press.
\bigskip\noindent
Kuo, H.-H. (1975).
{\it Gaussian measures in Banach Spaces.}
Lecture Notes in Mathematics 463.
Springer.
\bigskip\noindent
Landkof, N.S. (1972). {\it Foundations of Modern Potential Theory.} Springer.
\bigskip\noindent
Laha, R.G. and Rohatgi, V.K. (1979).
{\it Probability Theory}.
Wiley.
\bigskip\noindent
LeGall, J.-F. (1987). The exact Hausdorff measure of Brownian multiple points. In {\it Seminar on Stochastic Processes 1986} (E. Cinlar, K.L. Chung, R.K. Getoor eds.) Birkh\"auser.
\bigskip\noindent
Madrecki, A. (1983). On Gaussian type measures in p-adic Banach spaces. Preprint.
\bigskip\noindent
Madrecki, A. (1985). On Sazonov type topology in p-adic Banach space. {\it Math. Zeit.} {\bf 188} 223-236.
\bigskip\noindent
Madrecki, A. (1990)
Minlos' theorem in non-Archimedean locally compact convex spaces.
{\it Comment. Math. Prace Mat.} {\bf 30} 101-111.
\bigskip\noindent
Madrecki, A. (1991)
Some negative results on existence of Sazonov topology in $\ell$-adic Frechet spaces. {\it Arch. Math.} {\bf 56} 601-610.
\bigskip\noindent
Mahler, K. (1980). {\it p-adic Numbers and their Functions.} Cambridge University Press.
\bigskip\noindent
Marcus, M.B. and Pisier, G. (1984). {\it Random Fourier Series with Applications to Harmonic Analysis} Annals of Math. Studies 101. Princeton University Press.
\bigskip\noindent
Missarov, M.D. (1989). Random fields on the adele ring and Wilson's renormalization group. {\it Ann. Inst. Henri Poincar\'e Phys. Theor.} {\bf 50} 357-367.
\bigskip\noindent
Missarov, M.D. (1991). Renormalization group and renormalization theory in $p$-adic and adelic scalar models. In {\it Dynamical Systems and Statistical Mechanics (Moscow 1991), Adv. Soviet Math., 3.} American Mathematical Society.  
\bigskip\noindent
Monna, A. (1970). {\it Analyse Non-Archim\'edienne.} Springer.
\bigskip\noindent
Port, S.C. and Stone, C.J. (1978). {\it Brownian Motion and Classical
Potential Theory.} Academic Press.
\bigskip\noindent
Preston, C. (1971). Banach spaces arising from some integral inequalities. {\it Indiana Univ. Math. J.} {\bf 20} 997-1015.
\bigskip\noindent
Ruelle, Ph. and Thiran, E. (1989). Quantum mechanics on $p$-adic fields. {\it J. Math. Phys.} {\bf 30} 2854-2874.
\bigskip\noindent
Schikhof, W.H. (1984). {\it Ultrametric Calculus.} Cambridge University Press.
\bigskip\noindent
Sharpe, M. (1988). {\it General Theory of Markov Processes.} Academic Press.
\bigskip\noindent
Spokoiny, B.L. (1989). Non-Archimedean geometry and quantum mechanics. {\it Phys. Lett. B} {\bf 221} 120-124.
\bigskip\noindent
Taibleson, M.H. (1975). {\it Fourier Analysis on Local Fields.} Princeton University Press.
\bigskip\noindent
Taylor, S.J. (1973). Sample path properties of processes with stationary independent increments. In {\it Stochastic Analysis} (D.G. Kendall, E.F. Harding eds.) Wiley.
\bigskip\noindent
Taylor, S.J. and Wendel, J.G. (1966). The exact Hausdorff measure of the zero set of a stable process. {\it Z. Wahrscheinlichkeitstheorie verw. Geb.} {\bf 6} 170 - 180.
\bigskip\noindent
Trotter, H. (1958). A property of Brownian motion paths. {\it Illinois J. Math.} {\bf 2} 425-433.
\bigskip\noindent
van Rooij, A. (1978). {\it Non-Archimedean Functional Analysis.} Marcel Dekker.
\bigskip\noindent 
Vladimirov, V.S. and Volovich, I.V. (1989). p-adic quantum mechanics. {\it Commun. Math. Phys.} {\bf 123} 659-676.
\bigskip\noindent
Williams, D. (1979).
{\it Diffusions, Markov Processes and Martingales}.
Wiley.

\bye